\documentclass[twocolumn]{autart}    
\usepackage{amsmath} 
\usepackage{mathtools}
\usepackage{amssymb} 
\usepackage{upgreek}

\usepackage{graphicx}          

\usepackage{natbib}     
\usepackage{xcolor}
\usepackage{comment}

\usepackage{wrapfig}

\theoremstyle{definition}
\newtheorem{definition}{Definition}[section]

\begin{document}

\begin{frontmatter}

\title{Geometric sliding mode control of mechanical systems on Lie groups}

\thanks[footnoteinfo]{This paper was not presented at any IFAC 
meeting. $^{*}$Corresponding author Yu Tang, on leave from the National Autonomous University of Mexico, Mexico City, MEXICO.}

\author[UNAM]{Eduardo Espindola}\ead{eespindola@comunidad.unam.mx},
\author[UNAM]{Yu Tang$^{*}$}\ead{tang@unam.mx}               

\address[UNAM]{Ningbo Institute of Technology, Zhejiang University, Ningbo, CHINA.}  

\begin{keyword}                           
Geometric control; Lie groups; Mechanical systems; Sliding subgroups.               
\end{keyword}                             

\begin{abstract}                          

This paper presents a generalization of conventional sliding mode control designs for systems in Euclidean spaces to fully-actuated simple mechanical systems whose configuration space is a Lie group for the trajectory-tracking problem. 
A generic kinematic control is first devised in the underlying Lie algebra, which enables the construction of a Lie group on the tangent bundle where the system state evolves. A sliding subgroup is then proposed on the tangent bundle with the desired sliding properties, and a control law is designed for the error dynamics trajectories to reach the sliding subgroup globally exponentially. Tracking control is then composed of the reaching law and sliding mode, and is applied for attitude tracking on the special orthogonal group $SO(3)$ and the unit sphere $\mathcal S^3$. Numerical simulations show the performance of the proposed geometric sliding-mode controller (GSMC) in contrast with two control schemes of the literature.

\end{abstract}

\end{frontmatter}

\section{Introduction}
Sliding mode control (SMC) \citep{utkin1977variable} has been proven to be a very powerful control design method for systems evolving in Euclidean spaces. Its design usually consists of two stages: the reaching stage where the controller drives the system trajectories to a sliding surface, a subspace embedded in the Euclidean space designed to convoy some specific characteristics (e.g., convergence time, actuator saturation) in accordance with the given control objectives, and a sliding stage where the system trajectories converge to the origin according to the reduced-order dynamics constrained in the sliding surface, achieving the control objectives. In the sliding stage, the reduced-order dynamics is independent of the system dynamics, and therefore, this control design method ensures its robustness against a certain class of disturbances and has achieved great success in a wide range of applications. 

When this method is extended to mechanical systems whose configuration space is a general Lie group, care must be taken in the design of the sliding surface. Unlike the Euclidean case, when the system configuration space is a Lie group $G$, its time rate of change belongs to the tangent space $T_{g}G$ at the configuration $g$. Therefore, the state space is composed of the tangent bundle $G\times T_{g}G$. The topological structure and the underlying properties of the configuration space and the tangent space are very different. Without taking this into account in the SMC design, the sliding surface may not belong to the tangent bundle, and therefore no guarantee is offered to ensure that the system trajectories reach the sliding surface and the sliding mode may not exist at all \citep{cortes2019sliding}. The main problem is thus how to devise a group operation such that the tangent bundle is a Lie group and that the sliding subgroup is immersed in the tangent bundle so that the salient features of SMC in the Euclidean space mentioned above may be inherited by a general Lie group.

We present in this paper a general method of designing a sliding mode control, a geometric sliding mode control (GSMC),  for fully-actuated mechanical systems whose configuration space is a Lie group. A generic kinematic control is first devised in the underlying Lie algebra (the tangent space at the group identity with a bilinear map), which enables us to build a Lie group on the tangent bundle where the system state evolves. Then a sliding subgroup is proposed on the tangent bundle, and the sliding mode is guaranteed to exist. The sliding subgroup is designed to convoy control objectives, in particular, the almost global asymptotic convergence of the trajectories of the reduced-order dynamics to the identity of the tangent bundle is considered, which is the strongest convergence that may be achieved by continuous time-invariant feedback in a smooth Lie group \citep{bhat2000topological}. The reaching control law is then designed to drive the trajectories to the sliding subgroup globally exponentially. Tracking is then composed of the reaching law and sliding mode, as in the Euclidean case. 

\subsection{Related work}
The geometric approach to control designs has achieved significant advances for mechanical systems on nonlinear manifolds, for recent developments in this topic, see, for instance, \cite{bullo2019geometric} and the references therein. As recognized in \cite{koditschek1989autonomous,bullo1995control,maithripala2006almost}, a key point in control design is how to define the tracking error. The tracking error defined on a Riemannian manifold relying on an error function and a transport map  in \cite{bullo1999tracking} may be simplified if the manifold is endowed with a Lie group structure \citep{maithripala2006almost}, where the error notion can be globally defined explicitly  and is easier to be manipulated for stability analysis of the closed-loop system  
\citep{maithripala2015intrinsic,saccon2013optimal,de2018output,lee2012exponential,sarlette2010coordinated}. A similar situation is encountered in observer designs using an estimation error defined on a Riemannian manifold \citep{aghannan2003intrinsic} versus an estimation error defined by the group operation on Lie groups \citep{}{bonnabel2009non}.
The ability to define a global error on Lie groups provides a powerful tool for treating the error as an object in the state space globally and controlling it as a physical system so that the tracking problem can be reduced to stabilizing the error dynamics to the group identity \citep{bullo1995control,maithripala2006almost,spong2005controlled}.
 Moreover, a separation principle can be proved in the geometric approach to control designs \citep{maithripala2006almost,maithripala2015intrinsic} when part of the state in the control law is estimated by an exponentially convergent observer designed on the Lie group \citep{bonnabel2009non}, similar to an LTI system. This opens a wide field of applications for systems on Lie groups, such as rigid body motion control and trajectory tracking in $2$D and $3$D spaces, given the significant advances in both geometric control designs \citep{bullo2019geometric,spong2005controlled,lee2011geometric,akhtar2020controller,rodriguez2022new} and observer designs \citep{aghannan2003intrinsic,bonnabel2009non,mahony2008nonlinear,lageman2009gradient,zlotnik2018gradient}.

GSMC on Lie groups has been considered using two main approaches: developing the SMC in the underlying Lie algebra or developing it on the Lie group itself. The main idea in the former approach is 
first expressing the tracking error defined on the Lie group in its Lie algebra through the locally diffeomorphic logarithmic map \citep{bullo1995control}. Since the Lie algebra is a vector space,  a sliding surface can be designed as in the Euclidean case  \citep{culbertson2021decentralized,liang2021geometric,espindola2022novel}. 
In the latter approach, the sliding subgroup is designed directly on the Lie group. Since the topological structures of the configuration space (a Lie group) and the tangent space (a vector space) are very different when the underlying Lie group is not diffeomorphic to an Euclidean space, an important question arises as to how to ensure the sliding surface to be indeed a subgroup of the state space formed by the tangent bundle to guarantee the existence of the sliding mode and thus to inherit the salient features of SMC in the Euclidean space.

SMC designs using the second approach have been reported for Lie groups such as $SO(3)$, $\mathcal S^3$ for attitude control, and $SE(3)$ for motion controls \citep{ghasemi2020robust,lopez2021sliding}. However, the issue of whether the tangent bundle is a Lie group and whether the sliding subgroup is properly immersed on the tangent bundle was not addressed in these works. Therefore, the potential problem of lack of robustness due to the nonexistence of the sliding mode might appear. Recently,  \cite{cortes2019sliding} brought this issue to the attention of the control community, and proposed an SMC on the rotation group $SO(3)$ with a sliding surface which was ensured to be a Lie subgroup immersed in the tangent bundle $SO(3)\times \mathbb R^3$,  and a finite-time convergent controller was devised for attitude control. This design method was applied in \cite{meng2023second} to design a second-order SMC for fault-tolerant control designs.

\subsection{Contributions}
We generalize the conventional sliding mode control designs for systems in Euclidean spaces to fully-actuated simple mechanical systems whose configuration space is a Lie group for the trajectory-tracking problem. 
The main contributions can be summarized as follows: (1)  we endow the state space formed by the tangent bundle of the error dynamics with a Lie group structure by defining a group operation that is based on a generic kinematic control designed in the Lie algebra of the configuration Lie group;  (2) we design a smooth sliding subgroup and show it to be a Lie subgroup of the tangent bundle, therefore, inheriting the Lie group structure of the state space; and (3) we design a coordinate-free geometric sliding mode controller for a fully-actuated mechanical system on a Lie group which drives the error dynamics to the sliding subgroup globally exponential at the reaching stage, the error dynamics then converges  to the identity of the tangent bundle almost globally asymptotically at the sliding stage. In addition, rigid body tracking in $3$D space is addressed on the special orthogonal groups $SO(3)$ and on the unit sphere $\mathcal S^3$, respectively,  by applying the proposed geometric sliding mode control.


\subsection{Organization}
The rest of the paper is organized as follows. Section \ref{sec:background} presents the notation and background materials for simple mechanical systems with Lie groups as the configuration space. Section \ref{sec:TangentBundle-SlidingSurface} first endows the state space formed by the tangent bundle with a Lie group structure under a group operation, which is defined based on a generic kinematic control law in the Lie algebra of the configuration space. Then, a smooth sliding subgroup is defined, which is a Lie subgroup immersed in the tangent bundle. The convergence to the identity of the tangent bundle of the reduced-order dynamics constrained on the sliding surface is analyzed based on Lyapunov stability. Section \ref{GSMC} gives the design of the GSMC,  composed of a reaching law to the sliding subgroup and the convergence property of the sliding subgroup.  Attitude tracking of a rigid body in $3$D space is addressed  in Section \ref{AttitudeTracking} respectively on the rotational group $SO(3)$ and the unit sphere $\mathcal S^3$, and simulation results under the GSMC developed on $SO(3)$ are presented in Section \ref{sec:simulations} for illustration and comparison. Conclusions are drawn in Section \ref{sec:Conclusions}.

\section{Mechanical systems on Lie groups} \label{sec:background}

This section provides the notation and introduces the motion equations for a fully-actuated simple mechanical system on Lie groups. More details can be found in \cite{bullo2019geometric} and \cite{abraham2012manifolds}. 

Given a finite-dimension Lie group $G$, the identity of the group is denoted by $e\in G$. $T_{e}G$ denotes the tangent space in the identity, which also defines its Lie algebra $\mathfrak{g} \triangleq T_{e} G$ in the Lie bracket $[\cdot,\cdot]\in\mathfrak{g}$.  Let $L_{g}(h)=gh \in G$ and $R_{g}(h)=hg \in G$ be the left and right translation maps, respectively, $\forall g,h\in G$, and denote its corresponding tangent maps $ T_{e}L_{g}(\nu) = g \cdot \nu \in T_{g}G$ and $ T_{e}R_{g}(\nu) = \nu \cdot g \in T_{g}G$,  $\forall \nu \in T_{e}G$, it describes the natural isomorphism $T_{e}G \simeq T_{g}G$, which induces the equivalence $TG \simeq G\times T_{e}G$ for the tangent bundle $TG = G\times T_{g}G$. The inverse tangent map from $T_{g}G$ to $T_{e}G$ is denoted by $\nu = g^{-1}\cdot v^{L}_{g}$, where $v^{L}_{g} = \nu_{L}(g)\in T_{g}G$, being $\nu_{L}\in \Gamma^{\infty}(TG)$ a left-invariant vector field, with $\Gamma^{\infty}(TG)$ denoting the set of $C^{\infty}$-sections of $TG$, and respectively for a right-invariant vector field $\nu_{R}\in \Gamma^{\infty}(TG)$, it follows that $v^{R}_{g}=\nu_{R}(g) \in T_{g}G$, and accordingly $\nu = v^{R}_{g} \cdot g^{-1}$. 

The cotangent space at $g\in G$ is denoted by $T^{*}_{g}G$, while $\mathfrak{g}^{*}$ describes the dual space of the Lie algebra $\mathfrak{g}$. Likewise, the cotangent bundle is denoted by $T^{*}G \simeq G\times \mathfrak{g}^{*}$. Given a $\mathbb{R}$-vector space $V$, its dual space $V^{*}$, and a bilinear map $B\vcentcolon V\times V \to \mathbb{R}$, the flat map $B^{\flat}\vcentcolon V\to V^{*}$ is defined as $ \langle B^{\flat}(v);u \rangle =B(u,v)$, $\forall u,v \in V$, $B^{\flat}(v)\in V^{*}$, where $\langle \alpha ; v \rangle = \alpha(u)$ denotes the image in $\mathbb{R}$ of $v\in V$ under the covector $\alpha \in V^{*}$. If the flat map is invertible, then the inverse, known as the sharp map,  is denoted by $B^{\sharp}\vcentcolon V^{*} \to V$

The inner product on a smooth manifold $\mathcal{M}$ is denoted by $\langle\langle \cdot, \cdot\rangle\rangle \in \mathbb{R}$. A Riemannian metric $\mathbb{G}$ on a Lie group $G$ assigns the inner product $\mathbb{G}(g)\cdot \left( X_{g},Y_{g}\right)$ on each $T_{g}G$, $\forall X_{g},Y_{g}\in T_{g}G$. Moreover, when $\mathbb{G}$ is left-invariant (resp. right-invariant), it induces an inner product in the Lie algebra $\mathfrak{g}$ by $\mathbb{I}\left( \xi, \zeta \right) = \mathbb{G}(g)\cdot \left( \xi_{L}(g),\zeta_{L}(g)\right)$, $\forall \xi,\zeta \in \mathfrak{g}$. The kinetic energy is given by $\mathrm{KE}(v_{g})=(1/2)\mathbb{G}(g)\cdot \left( v_g,v_g\right) = (1/2) \mathbb{I} (\nu,\nu)$, where $\mathbb{I}$ is the kinetic energy tensor, which induces a kinetic energy metric $\mathbb{G}$ on $G$. In the rotational motion of a rigid body, $\mathbb{I}$ also represents the inertia tensor.

In the sequel, only the left invariance will be used. The proposed control methodology can be developed similarly for the right invariance. Also, subscripts and superscripts $L$ will be dropped when the meaning is clear. A left-invariant covariant derivative (affine connection) on a Lie group is denoted by $\nabla_{\xi_{L}}\zeta_{L}\in \Gamma^{\infty}(TG)$ for any vector fields $\xi_{L},\zeta_{L}\in \Gamma^{\infty}(TG)$. In addition, the Levi-Civita connection associated with the Riemannian metric $\mathbb{G}$ is denoted by $\overset{\mathbb{G}}{\nabla}$, 
which is unique and torsion-free. A left-invariant affine connection on a Lie group is uniquely determined by a bilinear map $B:\mathfrak{g}\times \mathfrak{g} \to \mathfrak{g}$ called the restriction of the left-invariant connection. In particular, the restriction for the left-invariant Levi-Civita connection  $\overset{\mathbb{G}}{\nabla}$ 
is defined as
\begin{equation}\label{eq:Rest}
    \overset{\mathfrak{g}}{\nabla}_{\xi}\zeta \triangleq \frac{1}{2}\left[ \xi,\zeta \right] - \frac{1}{2}\mathbb{I}^{\sharp}\left( \mathrm{ad}^{*}_{\xi}\mathbb{I}^{\flat}(\zeta) + \mathrm{ad}^{*}_{\zeta}\mathbb{I}^{\flat}(\xi) \right), 
\end{equation}
where  the adjoint map $\mathrm{ad}\vcentcolon \mathfrak{g}\times\mathfrak{g} \to \mathfrak{g}$ is defined as $\mathrm{ad}_{\xi}\zeta =\left[ \xi,\zeta \right]$, and $\mathrm{ad}^{*}_{\xi}\vcentcolon \mathfrak{g}^{*}\to \mathfrak{g}^{*}$ is the dual map defined as $\langle \mathrm{ad}^{*}_{\xi}\alpha ;\zeta  \rangle = \langle \alpha ; \left[\xi,\zeta\right]  \rangle$. 
Furthermore, the adjoint action $\mathrm{Ad}\vcentcolon G\times\mathfrak{g}\to \mathfrak{g}$ is  $\mathrm{Ad}_{g}\zeta = g\cdot \zeta \cdot g^{-1}$, $\forall g\in G$.
So, the left-invariant Levi-Civita connection is explicitly expressed as
\begin{equation}\label{eq:LeviC}
    \overset{\mathbb{G}}{\nabla}_{\xi_{L}}\zeta_{L} \triangleq \left( \mathrm{d}\zeta(\xi) +  \overset{\mathfrak{g}}{\nabla}_{\xi}\zeta \right)_{L},
\end{equation}
where $\mathrm{d}\zeta(\xi) \triangleq \frac{d}{dt}|_{t=0} \;\zeta\left( g\;\mathrm{exp}(\xi t)\right)$, 
being $\mathrm{exp}\vcentcolon \mathfrak{g}\to G$ the exponential map on $G$, which is a local $C^{\infty}$-diffeomorphism, and whose inverse is called the logarithmic map denoted by $\mathrm{log}\vcentcolon G\to \mathfrak{g}$. By the left-invariance of vector fields $\xi_{L},\zeta_{L}\in \Gamma^{\infty}(TG)$, the covariant derivative \eqref{eq:LeviC} is expressed in  terms of $\xi,\zeta\in\mathfrak{g}$ as follows
\begin{equation*}
    \nabla_{\xi}\zeta \triangleq \mathrm{d}\zeta(\xi) +  \overset{\mathfrak{g}}{\nabla}_{\xi}\zeta .
\end{equation*}


Consider a differentiable curve $g:  I\to G$, where $I$ is the set of all intervals. Then a body velocity $\nu:  I\to \mathfrak{g}$ is defined as $t\mapsto T_{g(t)}L_{g^{-1}(t)}\left( \dot{g}(t)\right)$, for all $t\in I$, and therefore
\begin{equation}\label{eq:Kin}
    \dot{g}(t) = g(t)\cdot \nu (t).
\end{equation}

A forced mechanical system is governed by the intrinsic Euler-Lagrange equations
\begin{equation}\label{eq:ELeq}
   \overset{\mathbb G}{\nabla}_{\dot{g}(t)}\dot{g}(t) = F_{u} + \Delta_{d}, 
\end{equation}
where $F_{u}=\sum^{m}_{a=1}u^{a}(t)\mathbb{G}^{\sharp}\left( T^{*}_{g(t)}L_{g^{-1}(t)}\left( f^{a}\right)\right)$ is the control force applied to the system on $T_{g}G$, being $u^{a}: I\to \mathbb{R}$ the control inputs, and $f^{a}(g)\in \mathfrak{g}^{*}$ the control forces. Furthermore, $\Delta_{d} \in T_{g}G$ represents the vector field version of constraint forces, such as potential external forces, uncontrolled conservative plus dissipative forces, and unmodeled disturbances. 

 In view of \eqref{eq:Kin} and the left-invariance of $\dot{g}$, the Levi-Civita connection in \eqref{eq:ELeq} can be  explicitly expressed using \eqref{eq:Rest}-\eqref{eq:LeviC} as
\begin{equation*}
   \overset{\mathbb G}{\nabla}_{\dot{g}(t)}\dot{g}(t)  = g(t)\cdot \left( \dot{\nu}(t) + \overset{\mathfrak{g}}{\nabla}_{\nu (t)}\nu (t)  \right), 
\end{equation*}
resulting  in the controlled Euler-Poincar\'e equation 
\begin{equation}\label{eq:EulerP}
    \dot{\nu}(t) + \overset{\mathfrak{g}}{\nabla}_{\nu (t)}\nu (t) = f_{u} + \delta_{d},
\end{equation}
with $f_{u}=\sum^{m}_{a=1}u^{a}(t)\mathbb{I}^{\sharp}\left(f^{a}\right) \in \mathfrak{g}$, and $\delta_{d}=g^{-1}\cdot \Delta_{d} \in \mathfrak{g}$.

The underlying mechanical system on the Lie group $G$ is then defined by the configuration Lie group $G$, the inertia tensor $\mathbb I$, and the external forces $f_u + \delta_d$.

\section{Lie Group Structure of the State Space and the Sliding Subgroup} \label{sec:TangentBundle-SlidingSurface}
In this section, we will endow the tangent bundle $TG\simeq G\times \mathfrak{g}$ with a Lie group structure by a properly designed group operation. For this purpose, an intrinsic control for kinematics is first proposed \eqref{eq:Kin}. Then we design a smooth sliding subgroup that is immersed in the tangent bundle so that it inherits the Lie group structure of the state space.  

\subsection{Intrinsic kinematic control}
The purpose of this subsection is to design a control law $\nu(t)\in\mathfrak{g}$ for the kinematics \eqref{eq:Kin} to render $g(t)\to e$, the group identity. Let $V:  G\to \mathbb{R}_{\geq 0}$ be an infinitely differentiable proper Morse function, which satisfies $V(g)>0$, $\forall g\in G\backslash \{ e\}$, $\mathrm{d}V(g)=0$ and $ V(g)=0 \iff g=e$. Morse functions, a class of error functions \citep{koditschek1989autonomous,bullo1995control}, are guaranteed to exist on many Lie groups of practical interest considered in this paper \citep{maithripala2015intrinsic,bullo2019geometric}. They represent potential energy that can be used to measure the distance between the configuration $g$ and the identity $e$ on $G$. The following definition specifies the class of kinematic controls considered in the paper.

\begin{definition}[Kinematic control law]\label{def:KctrLaw}
Let $g: I\to G$ be a differentiable curve governed by \eqref{eq:Kin}, for all $t\in I$. A kinematic control law is a map $\nu_{u}: G \to \mathfrak{g}$ that satisfies the following properties.
\begin{enumerate}\renewcommand{\theenumi}{\roman{enumi}}
    \item $\nu_{u}(e)=0$,\label{Pty1_DefKLaw}
    \item $\nu_{u}\left(g^{-1}\right) = -\nu_{u}\left(g\right)$,\label{Pty2_DefKLaw}
    \item $\langle \mathrm{d}V(g(t)); -g(t)\cdot \nu_{u}(g(t)) \rangle < 0$, $\forall g(t)\in G \backslash \mathcal{O}_{u}$, where $\mathcal{O}_{u}\triangleq \left\{ g\in G\backslash \{ e \} \; |\; g\cdot \nu_{u}(g) =0 \right\}$,\label{Pty3_DefKLaw}
    \item $\langle \mathrm{d}V(g(t)); -g(t)\cdot \nu_{u}(g(t)) \rangle = -y\left(g(t)\right)V\left(g(t)\right) $, $\forall g(t)\in\mathcal{U}$, where $y: \mathcal{U}\to \mathbb{R}_{>0}$, and $\mathcal{U}\subset G\backslash \mathcal{O}_{u}$ is a neighborhood of $e$. 
    \label{Pty4_DefKLaw}
\end{enumerate}
\end{definition}

Some comments on the class of kinematic controls are in order. Properties (\ref{Pty1_DefKLaw})-(\ref{Pty2_DefKLaw}) are instrumental to building a particular Lie subgroup on the tangent bundle. Properties (\ref{Pty3_DefKLaw})-(\ref{Pty4_DefKLaw}) represent the sliding (convergence) property of the reduced-order dynamics on the sliding subgroup (Lemma \ref{lem:SldLieSubG} below). In particular, Property (\ref{Pty3_DefKLaw}) states the almost-global asymptotic stability for system \eqref{eq:Kin} in closed loop with the kinematic control law $\nu (t)= -\nu_{u}(g(t))$. Note that since $\mathcal{O}_{u}$ is the set of closed-loop equilibria other than $g(t)=e$,  they are critical points of $V(g)$. Since $V(g)$ is a Morse function, the set $\mathcal{O}_{u}$ consists of a finite number of isolated points. In addition, this set is nowhere dense, which means that it cannot separate the configuration space. Therefore,  the complement $G\backslash \mathcal{O}_{u}$ is open and dense, i.e., $G\backslash \mathcal{O}_{u}$ is a submanifold of $G$ \citep{maithripala2006almost}. 
Finally, Property (\ref{Pty4_DefKLaw}) establishes the local exponential stability of the closed-loop system, where the existence of the neighborhood $\mathcal{U}$ is immediate because $V(g)$ is a Morse function, which has a unique minimum at $e\in G$ by definition.

Note that both $V(g)$ and $\nu_{u}(g)$ are of free design, provided that the properties in Definition \ref{def:KctrLaw} hold. However, it is worth considering the kinematic control law in the logarithmic coordinate, that is, $\nu_{u}(g)=\mathrm{log}(g)$, or some parallel vectors to $\mathrm{log}(g)$ \citep{akhtar2020controller}, as this map has been found to provide the strongest stability results, for example, almost global and local exponential convergence to the identity through a geodesic path \citep{bullo1995control}. 

\subsection{Lie Group structure for the state space}
For systems described in \eqref{eq:Kin} and \eqref{eq:EulerP}, the state space is the tangent bundle $TG \simeq G\times\mathfrak{g}$. To endow it with a Lie group structure, we consider the binary operation $\star:\ TG\times TG\mapsto TG$  defined in the following
\begin{align}\label{eq:BnyOp}
    &h_{1}\star h_{2} \triangleq \notag\\ 
    &\big( g_{1}g_{2} ,\; \nu_{1} + \nu_{2} + \lambda\nu_{u}(g_{1}) +\lambda \nu_{u}(g_{2}) -\lambda \nu_{u}(g_{1}g_{2}) \big), 
\end{align}
$\forall h_{1}=(g_{1},\nu_{1}),\; h_{2}=(g_{2},\nu_{2}) \in TG$, and $\lambda\in \mathbb R_{>0}$.

\begin{lem}[The state space $TG$ as a Lie group]\label{lem1:LieGroup}
The tangent bundle $TG \equiv G\times\mathfrak{g}$ endowed with the binary operation \eqref{eq:BnyOp} is a Lie group, with
\begin{enumerate}\renewcommand{\theenumi}{\roman{enumi}}
    \item Identity element: $f \triangleq (e,0)\in TG$,\label{Pty1_ThmLG}
    \item Inverse element:  $h^{-1}\triangleq\left( g^{-1},-\nu\right)\in TG$, $\forall h=(g,\nu) \in TG$. \label{Pty2_ThmLG}
\end{enumerate}
\end{lem}
\begin{pf}
Being $TG$ a smooth manifold with \eqref{eq:BnyOp} a smooth operation, it only remains to verify the group axioms as follows.
\begin{enumerate}
    \item $\forall h=(g,\nu)\in TG$, it satisfies
    \begin{align*}
        h\star f &=  \left( ge,\; \nu + 0 + \lambda\nu_{u}(g) +\lambda \nu_{u}(e) -\lambda \nu_{u}(ge)  \right) \notag\\
        &= f\star h = h,
    \end{align*}
    where Property  of Definition \ref{def:KctrLaw}(\ref{Pty1_DefKLaw})  is used.
    \item The group operation between $h=(g,\nu)\in TG$ and its inverse  $h^{-1}=\left( g^{-1},-\nu\right)\in TG$  verifies 
    \begin{align*}
        &h^{-1}\star h\\
        &= \big( g^{-1}g, -\nu + \nu + \lambda\nu_{u}(g^{-1}) +\lambda \nu_{u}(g) -\lambda \nu_{u}(g^{-1}g) \big) \notag\\
        &=\big( gg^{-1}, \nu - \nu + \lambda \nu_{u}(g) +\lambda\nu_{u}(g^{-1}) 
         -\lambda \nu_{u}(gg^{-1}) \big) \notag\\
        &= h \star h^{-1} = f,
    \end{align*}
    where Properties (\ref{Pty1_DefKLaw})-(\ref{Pty2_DefKLaw}) of Definition \ref{def:KctrLaw} are used.
    \item The associativity  $h_{1}\star \left( h_{2}\star h_{3}\right) = \left( h_{1}\star h_{2}\right) \star h_{3}$ is proved straightforwardly by substitution, using the properties of Definition \ref{def:KctrLaw}. 
\end{enumerate}
\hfill$\blacksquare$
\end{pf}

\begin{rem}[Tangent bundle $TG$]\label{rmk:Lem1}
The definition of the group operation  \eqref{eq:BnyOp} relying on the kinematic control $\nu_{u}(g)$ in Definition \ref{def:KctrLaw} is crucial to define a sliding Lie subgroup immersed  in $TG$ in the next subsection. In fact, a group operation to endow $TG$ to be a Lie group may simply be $h_{1}\star h_{2} = \left( g_{1}g_{2},\; \nu_{1}+\nu_{2}\right)$.  However, this operation does not allow to design of a useful sliding subgroup, in particular, it fails to prove closure under the group operation, as will be seen below.
\end{rem}

\begin{rem}[Associativity] 
The associativity proved in Lemma \ref{lem1:LieGroup} ensures the proposed Lie group $TG$ to be globalizable \citep{olver1996non}, that is, the local Lie group $TG$ can be extended to be a global topological group. This fact allows us to develop a sliding mode control defined globally on the state space in contrast to the Lie groups defined locally  in \cite{cortes2019sliding} and \cite{meng2023second}.

\end{rem}

\subsection{Sliding Subgroup on $TG$ }
In this subsection, we define a smooth sliding subgroup on the tangent bundle.
The following lemma shows that $H\subset TG$ is an immersed submanifold of $TG$ that inherits the topology and smooth structure of the tangent bundle $TG$ \citep{lee2013smooth}.

\begin{lem}[Sliding Lie subgroup]\label{lem2:LieSubGroup}
Define  
\begin{equation} 
    H\triangleq \left\{ h=(g,\nu)\in TG \;|\; s(h)=0\right\}\subset TG, \label{eq:slidingsurfaceH}
\end{equation}
 where  $\forall h=(g,\nu)\in TG$, the map $s: TG\mapsto \mathfrak{g}$ is defined as
\begin{equation}\label{eq:sldS}
    s(h) = \nu + \lambda \nu_{u}(g).
\end{equation}
Then $H\subset TG$ is a Lie subgroup under the group operation \eqref{eq:BnyOp}.
\end{lem}
\begin{pf}
The smoothness of $H$ is immediate, because the map defined in \eqref{eq:sldS} is smooth. The proof consists thus in showing that subset $H$ inherits the group structure of the Lie group  $TG$, by verifying the following: 
\begin{enumerate}\renewcommand{\theenumi}{\roman{enumi}}
    \item {\it Identity}:  The identity of the tangent bundle  $f=(e,0)\in H$. This is immediate by  Definition \ref{def:KctrLaw}(\ref{Pty1_DefKLaw}) since $s(f)=0+\lambda \nu_{u}(e) = 0$.
    \item {\it Inverse}. $\forall h=(g,\nu)\in H$, $s(h)=0\implies \nu = -\lambda\nu_{u}(g)$. By Definition \ref{def:KctrLaw}(\ref{Pty2_DefKLaw}) it follows that
    \begin{align*}
        s\left( h^{-1}\right) &= -\nu +\lambda \nu_{u}\left( g^{-1}\right) \notag\\
        &= -\left( -\lambda\nu_{u}(g) \right) +\lambda \nu_{u}\left( g^{-1}\right) =0.
    \end{align*}
    This proves that $h^{-1}\in H$ for all $h\in H$.
    \item {\it Closure}. Given  $h_{1}=(g_{1},\nu_{1}), \ h_{2}=(g_{2},\nu_{2})\in H$, then $s(h_{1})=0\implies \nu_{1}=-\lambda\nu_{u}(g_{1})$ and $s(h_{2})=0\implies \nu_{2}=-\lambda\nu_{u}(g_{2})$. By \eqref{eq:BnyOp} $h_1\star h_2=\big(g_1g_2, -\lambda\nu_u(g_1g_2)\big)$. Thus,  $s\left( h_{1}\star h_{2}\right)=-\lambda\nu_u(g_1g_2)+\lambda\nu_u(g_1g_2)=0$. That is, $H$ is closed under the group operation.
\end{enumerate}
\hfill$\blacksquare$
\end{pf}

The following lemma shows that once a trajectory reaches the sliding subgroup it will stay on it and converges to the group identity.

\begin{lem}[Properties of the sliding subgroup $H$] \label{lem:SldLieSubG}
Consider the sliding Lie subgroup $H\subset TG$ in \eqref{eq:slidingsurfaceH}. Then $H$ is forward 
invariant, i.e., $h(t_r)\in H$ for some $t_r\in I$ $\implies h(t)\in H, \forall t\geq t_r$. Moreover, $h(t)\to (e,\ 0)$ almost globally asymptotically.
\end{lem}
\begin{pf}
Consider a differentiable curve $g: I \to G$ of the dynamics \eqref{eq:Kin}.
Let $V:  G\to \mathbb{R}$ be a proper Morse function with the unique minimum  at $e\in G$. Then, along the trajectory $g(t)$ and  $\forall t\in I$, it yields
\begin{equation*}
    \frac{\mathrm{d}}{\mathrm{d}t}V\left(g(t)\right) = \langle \mathrm{d}V\left(g(t)\right); \dot{g}(t) \rangle = \langle \mathrm{d}V\left(g(t)\right); g(t)\cdot \nu(t) \rangle.
\end{equation*}

Assume that $h(t_r)\in H$, for some $t_r\in I$. Then $s(h)=0$ gives $\nu(t) = -\lambda\nu_{u}(g(t))$. Therefore,
\begin{align*}
    \frac{\mathrm{d}}{\mathrm{d}t}V\left(g(t)\right) &= \langle \mathrm{d}V\left(g(t)\right); g(t)\cdot \nu(t) \rangle \notag\\
    &= \langle \mathrm{d}V\left(g(t)\right); -\lambda g(t)\cdot \nu_{u}\left(g(t)\right) \rangle ,
\end{align*}
In light of Definition \ref{def:KctrLaw}(\ref{Pty3_DefKLaw}), it follows that $\frac{\mathrm{d}}{\mathrm{d}t}V\left(g(t)\right) <0$, for all $g(t_r)\in G\backslash \mathcal{O}_{u}$, and $\frac{\mathrm{d}}{\mathrm{d}t}V\left(g(t)\right) = 0 \iff g=e$, where $\mathcal{O}_{u}$ is a nowhere-dense set with a finite number of points given in Definition \ref{def:KctrLaw}(\ref{Pty3_DefKLaw}). Therefore, $h(t)$ will remain on $H$ for all $t\geq t_r$, and the equilibrium $g(t)=e$ of \eqref{eq:Kin} is almost globally asymptotically stable for all $g(t_r)\in G\backslash \mathcal{O}_{u}$ and locally exponentially stable $\forall g(0)\in \mathcal{U}$, according to Definition \ref{def:KctrLaw}(\ref{Pty4_DefKLaw}). 

\hfill$\blacksquare$
\end{pf}


\section{Geometric Sliding Mode Control (GSMC)} \label{GSMC}
  In this section, we design a control law, called the reaching law,  
  for $f_u$ in the Euler-Poincar\'e equation \eqref{eq:EulerP}
  to drive the trajectory $h(t)=(g(t),\nu (t))\in TG$ to the sliding subgroup $H$. 
  Then the tracking control objective will be achieved as a consequence of Lemma \ref{lem:SldLieSubG}.
  
\subsection{Reaching Law}

The Euler-Lagrange dynamics \eqref{eq:ELeq}, ignoring disturbance forces $\Delta_{d}$,  is expressed as
\begin{equation}\label{eq:ELintrinsic}
    \nabla_{\nu (t)}\nu(t) = f_{u}.
\end{equation}
which is defined on $TG$, being $h(t)=\left( g(t),\nu (t)\right)$ the state variable.

 The intrinsic acceleration for the sliding variable \eqref{eq:sldS} is calculated, by using \eqref{eq:ELintrinsic}, through the covariant derivative of $s(h)\in\mathfrak{g}$ with respect to itself as 
\begin{align*}
    \nabla_{s(h)}s(h) &= \frac{d}{dt}s(h) + \overset{\mathfrak{g}}{\nabla}_{s(h)}s(h) \notag \\
    &= \dot{\nu} + \lambda \dot{\nu}_{u}(g)+ \overset{\mathfrak{g}}{\nabla}_{s(h)}s(h).
\end{align*}
Substituting the Euler-Poincar\'e equation \eqref{eq:EulerP} yields
\begin{align}\label{eq:sp2}
    &\nabla_{s(h)}s(h)\notag\\
    =& -\overset{\mathfrak{g}}{\nabla}_{\nu (t)}\nu (t) + \lambda \dot{\nu}_{u}(g)+ \overset{\mathfrak{g}}{\nabla}_{s(h)}s(h) + f_{u} \notag\\
    =& \mathbb{I}^{\sharp}\left( \mathrm{ad}^{*}_{\nu (t)}\mathbb{I}^{\flat}(\nu (t)) \right) - \mathbb{I}^{\sharp}\left( \mathrm{ad}^{*}_{s(h)}\mathbb{I}^{\flat}(s(h)) \right) 
    + \lambda \dot{\nu}_{u}(g) + f_{u},
\end{align}
where  the skew-symmetry of the Lie bracket $[\cdot , \cdot]\in\mathfrak{g}$  in \eqref{eq:Rest} is used.  The reaching law is then proposed as follows
\begin{equation}\label{eq:Ctrl2}
    f_{u} = \mathbb{I}^{\sharp}\left( \mathrm{ad}^{*}_{\lambda\nu_{u}(g)}\mathbb{I}^{\flat}(\nu (t)) \right) -\lambda\dot\nu_{u}(g) -k_{s}s(h),
\end{equation}
with $k_{s}>0$ a design parameter.

\begin{thm}[Reaching Controller]\label{thm:Ctrl2}
  The reaching law \eqref{eq:Ctrl2} drives exponentially the trajectories of the closed-loop  system \eqref{eq:sp2} to the sliding subgroup $H$ $\forall h(0)\in TG$, i.e., $s(h(t))\to 0$ globally exponentially. 
\end{thm}
\begin{pf}
Consider the function $W: TG \to \mathbb{R}$ defined below
\begin{equation}\label{eq:V2}
    W(h) =\frac{1}{2}\mathbb{I}(s(h),s(h)).
\end{equation}

Its time evolution along trajectories of \eqref{eq:sp2} is given by
\begin{align*}
    \dot{W}(h) &= \mathbb{I}\left( \nabla_{s(h)}s(h),\; s(h)\right) \notag\\
    &= \mathbb{I}\left( \mathbb{I}^{\sharp}\left( \mathrm{ad}^{*}_{\nu (t)}\mathbb{I}^{\flat}(\nu (t)) \right) - \mathbb{I}^{\sharp}\left( \mathrm{ad}^{*}_{s(h)}\mathbb{I}^{\flat}(s(h)) \right) \right. \notag\\
    &\quad\quad + \lambda \dot{\nu}_{u}(g) + f_{u},\; s(h)\Bigl),
\end{align*}
which in closed loop with the controller \eqref{eq:Ctrl2} yields
\begin{align*}
    \dot{W}(h) &= \mathbb{I}\left( \mathbb{I}^{\sharp}\left( \mathrm{ad}^{*}_{\nu (t)}\mathbb{I}^{\flat}(\nu (t)) \right) - \mathbb{I}^{\sharp}\left( \mathrm{ad}^{*}_{s(h)}\mathbb{I}^{\flat}(s(h)) \right) \right. \notag\\
    &\quad\quad +\mathbb{I}^{\sharp}\left( \mathrm{ad}^{*}_{\lambda\nu_{u}(g)}\mathbb{I}^{\flat}(\nu (t)) \right)  -k_{s}s(h),\; s(h)\Bigl) ,\notag\\
    &=\mathbb{I}\left( \mathbb{I}^{\sharp}\left( \mathrm{ad}^{*}_{s(h)}\mathbb{I}^{\flat}(\nu (t)) \right) - \mathbb{I}^{\sharp}\left( \mathrm{ad}^{*}_{s(h)}\mathbb{I}^{\flat}(s(h)) \right) \right. \notag\\
    &\quad\quad -k_{s}s(h),\; s(h)\Bigl). 
\end{align*}
By Lemma \ref{lm:skysym} in Appendix \ref{AppA} the term $\mathbb{I}\left( \mathbb{I}^{\sharp}\left(\mathrm{ad}^{*}_{\zeta}\mathbb{I}^{\flat}(\eta)\right),\;\zeta\right) = 0$, for any $\zeta,\eta\in \mathfrak{g}$. Therefore,  
\begin{equation*}
    \dot{W}(h) = -k_{s}\mathbb{I}\left(s(h),\; s(h)\right)=-2k_s W(h).
\end{equation*}
It follows from Proposition 6.26 of \cite{bullo2019geometric} that $W(h(t))\to 0$ exponentially.
\hfill$\blacksquare$
\end{pf}

\begin{rem}[Passivity of the Lagrangian dynamics]\label{rmk:Ctrl1}
Note that the first two right-hand terms of the control law \eqref{eq:Ctrl2} complete the 
terms $\mathbb{I}^{\sharp}\mathrm{ad}^{*}_{s(h)}\mathbb{I}^{\flat}(\nu (t)) - \mathbb{I}^{\sharp} \mathrm{ad}^{*}_{s(h)}\mathbb{I}^{\flat}(s(h))$. By exploring the intrinsic passivity properties in Lemma \ref{lm:skysym} in Appendix \ref{AppA}, these terms were not canceled in the above stability analysis. This result was first given for the Lie group $SO(3)$ in \cite{koditschek1989autonomous}.  The lemma \ref{lm:skysym} extends this result to coordinate-free Lagrangian dynamics on a general Lie group, which has not been explored, to the authors' knowledge, in the literature for stability analysis.
\end{rem}

\begin{rem}[The reaching controller ]\label{rmk:Ctrl2}
The reaching law \eqref{eq:Ctrl2} achieves the convergence of $s(h(t))\to 0$ for the Euler-Lagrange dynamics \eqref{eq:ELintrinsic}, which implies that $h(t)\in TG$ reaches the 
sliding subgroup $H$ exponentially.  Note that the result of Theorem \ref{thm:Ctrl2} holds when the external constraint forces $\delta_{d}$ can be compensated for by the controller $f_{u}$, which was omitted from the control design. Otherwise, in the presence of bounded $\delta_{d}$, $h(t)\in TG$ will remain bounded and close to $H$. 
\end{rem}

\subsection{Tracking Control}

Let $g_{r}: I\to G$ be a twice differentiable configuration reference, with the corresponding reference body velocity $\nu_{r}: I\to \mathfrak{g}$ given by $\nu_{r} (t)\triangleq g^{-1}_{r}(t)\cdot \dot{g}_{r}(t)$. The problem is to design a control law  $f_u$ to track the reference.
The Lie group structure of the configuration space $G$ enables to define the following intrinsic configuration error 
\begin{equation*}
    g_{e}(t) \triangleq g^{-1}_{r}(t)g(t).
\end{equation*}
By left invariance  
the body velocity error is defined as  
\begin{align}\label{eq:velErr}
    \nu_{e}(t) &\triangleq g^{-1}_{e}(t)\cdot \dot{g}_{e}(t) = \nu (t) - \eta_{r}(t),
\end{align}
with $\eta_{r}(t)= \mathrm{Ad}_{g^{-1}_{e}} \nu_{r}(t)$. Then, the error dynamics evolving on $TG$ is described by 
\begin{equation}\label{eq:ELintrinsic-Nu-e}
    \nabla_{\nu_e (t)}\nu_e(t) = f_{u},
\end{equation}
being  the state variable $h_e(t)=\left( g_e(t),\nu_e (t)\right) \in TG$.

The tracking problem, therefore, boils down to stabilizing the identity $f=(e, 0)$ on $TG$. 
By using the sliding-model control strategy, the error state is first driven to the sliding subgroup in the reaching stage, and then on the sliding subgroup, the reduced-order dynamics converges to the identity $f$ ensured by Lemma \ref{lem:SldLieSubG}. 


In terms of the error state $h_e$  the sliding variable \eqref{eq:sldS} is given by
\begin{equation}\label{eq:SurfaceErr}
    s(h_{e}) = \nu_{e}(t) + \lambda \nu_{u}(g_{e}),
\end{equation}
and, its covariant derivative, by using \eqref{eq:Rest}-\eqref{eq:LeviC}, 
is
\begin{align*}  &\nabla_{s(h_{e})}s(h_{e})\\
    &= \frac{d}{dt}s(h_{e}) + \overset{\mathfrak{g}}{\nabla}_{s(h_{e})}s(h_{e}) \notag\\
    &= \dot{\nu}_{e}(t) + \lambda\dot{\nu}_{u}(g_{e}) -\mathbb{I}^{\sharp}\left( \mathrm{ad}^{*}_{s(h_{e})}\mathbb{I}^{\flat}\left( s(h_{e}) \right) \right) \notag\\
    &= \dot{\nu}(t) - \dot{\eta}_{r}(t) + \lambda\dot{\nu}_{u}(g_{e})
    -\mathbb{I}^{\sharp}\left( \mathrm{ad}^{*}_{s(h_{e})}\mathbb{I}^{\flat}\left( s(h_{e}) \right) \right).
\end{align*}
Ignoring disturbance $\delta_{d}$  it follows from  the Euler-Poincar\'e equation \eqref{eq:EulerP} that
\begin{align}
    \nabla_{s(h_{e})}s(h_{e})=&\mathbb{I}^{\sharp}\left( \mathrm{ad}^{*}_{\nu (t)}\mathbb{I}^{\flat}\left( \nu (t) \right) \right) + f_{u} - \dot{\eta}_{r}(t) \label{eq:ELerrDyn}\\
    &+ \lambda\dot{\nu}_{u}(g_{e}) -\mathbb{I}^{\sharp}\left( \mathrm{ad}^{*}_{s(h_{e})}\mathbb{I}^{\flat}\left( s(h_{e}) \right) \right) . \notag
\end{align}
We proposed the following tracking controller 
\begin{align}\label{eq:TckCtrl}
   f_{u} &=  \mathbb{I}^{\sharp}\left( \mathrm{ad}^{*}_{\lambda\nu_{u}(g_{e})-\eta_{r}(t)}\mathbb{I}^{\flat}(\nu (t)) \right) -\lambda\dot\nu_{u}(g_{e}) + \dot{\eta}_{r}(t) \notag\\
   &\quad -k_{s}s(h_{e}) ,
\end{align}
where $k_{s}>0$ is a design parameter. The following theorem establishes the stability of the equilibrium $h_e=f$ in the closed-loop system \eqref{eq:ELerrDyn}-\eqref{eq:TckCtrl}.

\begin{thm}[Tracking Controller]\label{thm:TckCtrl}
Consider the error dynamics \eqref{eq:ELerrDyn} in closed loop with the controller \eqref{eq:TckCtrl}. Then, the equilibrium  $h_{e}(t)=f$ is
\begin{enumerate}\renewcommand{\theenumi}{\roman{enumi}}
\item almost-globally asymptotically stable, for all $h_{e}(0)\in \overline{TG}\triangleq G\backslash\mathcal{O}_{u}\times \mathfrak{g}$,
\item locally exponentially stable for all $h_{e}(0)\in \overline{TU}\triangleq\mathcal{U}\times \mathfrak{g}$, where $\mathcal{O}_{u}$ and $\mathcal{U}$ are given in Definition \ref{def:KctrLaw}(\ref{Pty3_DefKLaw})-(\ref{Pty4_DefKLaw}).
\end{enumerate}
\end{thm}
\begin{pf}
Substituting the controller \eqref{eq:TckCtrl} in the error dynamics  \eqref{eq:ELerrDyn} yields  the closed-loop dynamics
\begin{align*}
    \nabla_{s(h_{e})}s(h_{e}) &= \mathbb{I}^{\sharp}\left( \mathrm{ad}^{*}_{s(h_{e})}\mathbb{I}^{\flat}\left( \nu (t) \right) \right) -k_{s}s(h_{e}) \notag\\
    &\quad  -\mathbb{I}^{\sharp}\left( \mathrm{ad}^{*}_{s(h_{e})}\mathbb{I}^{\flat}\left( s(h_{e}) \right) \right) ,
\end{align*}
which has an equilibrium point at $s(h_{e})=0$. The results follow as a consequence of Theorem \ref{thm:Ctrl2} (the reaching stage) and Lemma \ref{lem:SldLieSubG} (the sliding mode).
\hfill$\blacksquare$
\end{pf}

\begin{rem}[The tracking controller]\label{rmk:TckCtrl}

Theorem \ref{thm:TckCtrl} gives a coordinate-free sliding mode control for a mechanical system whose configuration space is a general Lie group. The group structure allows defining globally a tracking error, whose dynamics evolves on the tangent bundle. The Lie subgroup of the sliding subgroup immersed  on the tangent bundle ensures the existence of the sliding mode and thus inherits the salient features of the SMC in Euclidean spaces.

Similarly to the Euclidean case, the design of the sliding subgroup and the reaching law may incorporate other control objectives, such as finite-time convergence and controller saturation, which are, however, beyond the scope of the main purposes of this paper. 
\end{rem}

\section{Attitude Tracking of a Rigid Body}\label{AttitudeTracking}

In this section, we present the attitude tracking of a rigid body in the $3$D space using the proposed GSMC. To illustrate the theoretic development, the problem is addressed using attitude representation by first the rotation matrix on $SO(3)$ and then by the unit quaternion $\mathcal S^3$. 

\subsection{GSMC for Attitude Tracking on $SO(3)$}

The group of rotations on $\mathbb{R}^{3}$ is the Lie group $SO(3)=\left\{ R\in\mathbb{R}^{3\times 3}\; |\; RR^{T}=R^{T}R= I_{3}, \; \mathrm{det}(R)=+1\right\}$, with the usual multiplication of matrices as the group operation. The identity of the group is the identity matrix $I_{3}$ of $3\times 3$,  and the inverse is the transpose $R^{T}\in SO(3)$ for any $R\in SO(3)$. The Lie algebra is given by the set of skew-symmetric matrices $\mathfrak{so}(3)=\left\{ S\in\mathbb{R}^{3\times 3}\; |\; S^{T}=-S \right\}$, which is isomorphic to $\mathbb R^3$, i.e., $\mathfrak{so}(3)\simeq \mathbb{R}^{3}$. The Lie bracket in $\mathbb R^3$ is defined by the cross product $[ \zeta , \eta] = \mathrm{ad}_{\zeta}\eta\triangleq \zeta\times\eta$, $\forall \zeta,\eta\in\mathbb{R}^{3}$. Denote the isomorphism $\cdot^{\wedge}\vcentcolon \mathbb{R}^{3}\to \mathfrak{so}(3)$, and respectively the inverse map $\cdot^{\vee}\vcentcolon \mathfrak{so}(3) \to \mathbb{R}^{3}$. Then for a differentiable curve $R\vcentcolon I \to SO(3)$ with left-invariant dynamics $\dot{R}(t)\in T_{R}SO(3)$, the body angular velocity is given by  
\begin{equation*}
    \Omega^{\wedge}(t) = R^{T}(t)\dot{R}(t) = \left[ \begin{array}{ccc}
         0 & -\Omega_{3}(t) & \Omega_{2}(t)  \\
         \Omega_{3}(t) & 0 & -\Omega_{1}(t) \\ 
         -\Omega_{2}(t) & \Omega_{1}(t) & 0
    \end{array} \right] ,
\end{equation*}
for all $t\in I$. The kinetic energy of the rotational motion of a rigid body is calculated as $\mathrm{KE}(\Omega) = \frac{1}{2}\mathbb{J}(\Omega, \Omega) \triangleq   \frac{1}{2} \langle\langle\mathbb{J}\Omega,\Omega\rangle\rangle$, where $\mathbb{J}=\mathbb{J}^{T}\in\mathbb{R}^{3\times 3}$ is the positive-definite inertia tensor. Therefore, $\mathrm{ad}^{*}_{\zeta}\mathbb{J}^{\flat}(\eta) = \left(\mathbb{J}\eta\right)^{\wedge} \zeta$, and $\mathbb{J}^{\sharp}(\zeta) = \mathbb{J}^{-1}\zeta$. Hence, the rotational motion described by the Euler-Lagrange equation \eqref{eq:ELeq} is 
\begin{equation}\label{eq:ELDynSO3}
    \nabla_{\Omega(t)}\Omega(t) = \tau_{u}.
\end{equation}
The state $(R,\omega)$ evolves on the tangent bundle $TSO(3)\simeq SO(3)\times \mathbb{R}^{3}$, 
and the control torque $\tau_{u} = \mathbb{J}^{-1}\tau \in\mathbb{R}^{3}$ is  expressed in the body frame. Furthermore, \eqref{eq:ELDynSO3} is explicitly expressed, by using the Euler-Poincar\'e equation \eqref{eq:EulerP} and restriction \eqref{eq:Rest}, as
\begin{equation}\label{eq:DynSO3}
    \dot{\Omega}(t) - \mathbb{J}^{-1}\left(\mathbb{J}\Omega(t)\right)^{\wedge} \Omega(t) = \tau_{u}.
\end{equation}

Let $R_{r}\vcentcolon I\to SO(3)$ be a twice differentiable attitude reference, and $\Omega_{r}: I\to \mathbb{R}^{3}$,  the reference angular velocity expressed in the body frame, which holds $\Omega_{r}(t)=\left(R^{T}_{r}(t)\dot{R}_{r}(t)\right)^{\vee}$. 
Then, the intrinsic attitude error is
\begin{equation*}
    R_{e}(t) \triangleq R^{T}_{r}(t)R(t). 
\end{equation*}
In view of \eqref{eq:velErr} the (left-invariant) velocity error is 
\begin{align}
    \Omega_{e}(t) &\triangleq \left( R^{T}_{e}(t)\dot{R}_{e}(t) \right)^{\vee} = \Omega(t) - \sigma(t),\label{eq:ErrKynSO3}\\
    \sigma (t) &\triangleq \mathrm{Ad}_{R^{-1}_{e}}\Omega_{r}(t) = R^{T}_{e}(t)\Omega_{r}(t) \notag.
\end{align}
Therefore, the distance between $R_{e}(t)$ and $I_{3}$ is properly measured with the Morse function $V_{1}(R_{e})\triangleq 2-\sqrt{1+\mathrm{tr}(R_{e}(t))}$, proposed by \cite{lee2012exponential}. In fact, $V_{1}(R_{e})=0\iff R_{e}=I_{3}$ and is positive for all $R_{e}\in SO(3)\backslash \{ I_{3}\}$. Moreover, along the trajectories of \eqref{eq:ErrKynSO3}, it satisfies
\begin{align*}
    \frac{\mathrm{d}}{\mathrm{d}t}V_{1}(R_{e}) &= \left\langle\left\langle \psi(R_{e})\left(R_{e}(t)-R^{T}_{e}(t)\right)^{\vee}, \;\Omega_{e}(t)\right\rangle\right\rangle ,\notag \\
    \psi(R_{e}) &\triangleq \frac{1}{2\sqrt{1+\mathrm{tr}(R_{e})}},
\end{align*}
for all $R_{e}\in SO(3)\backslash \mathcal{O}_{R}$, where $\mathcal{O}_{R}\triangleq \{ R\in SO(3) | \mathrm{tr}(R)$ $=-1 \}$.
Furthermore, given $\mathcal{U}_{R}$ $\triangleq$ $\{ R_{e}\in SO(3)\backslash$ $\mathcal{O}_{R} \;|\; V_{1}(R_{e})< 2-\epsilon \}$, for some $\epsilon>0$ arbitrarily small, $V_{1}(R_{e})$ verifies \citep{lee2012exponential}
\vspace{-10mm}{\small
\begin{align*}
\left\|\psi(R_{e})\left(R_{e}-R^{T}_{e}\right)^{\vee}\right\|^{2} &\leq V_{1}(R_{e}) \leq 2\left\|\psi(R_{e})\left(R_{e}-R^{T}_{e}\right)^{\vee}\right\|^{2},
\end{align*}}
for all $R_{e}\in \mathcal{U}_{R}$.

Consider the kinematic control law 
\begin{align}
 \Omega_{u}(R_{e}) & \equiv \mathrm{log}(R_{e})^{\vee} ,\label{eq:KinCtrlLawSO3}\\
 \mathrm{log}(R_{e}) &\triangleq  
 \left\{ 
 \begin{array}{ll}
    0_{3\times 3},  & R_{e} = I_{3},   \\
    \frac{\phi(R_{e})}{2\sin\left( \phi(R_{e})\right)} \left( R_{e} - R^{T}_{e}\right) ,  & R_{e} \neq I_{3},
 \end{array} \right.  \notag
\end{align}
where $\phi(R_{e}) \triangleq \arccos{\left( \frac{1}{2}\left( \mathrm{tr}(R_{e})-1\right)\right)} \in (-\pi,\pi)$,
and $0_{n\times m}$ is a matrix of size $n\times m$ with zero-entries. It can verify readily Definition \ref{def:KctrLaw}(\ref{Pty1_DefKLaw})-(\ref{Pty2_DefKLaw})  by  \eqref{eq:KinCtrlLawSO3}. To verify 
Definition \ref{def:KctrLaw}(\ref{Pty3_DefKLaw})-(\ref{Pty4_DefKLaw}) under the kinematic control $\Omega_{e} (t)=-\Omega_{u}(R_{e})$, consider the derivative of the Morse function $V_{1}(R_{e})$ along the error kinematics $\dot R_e=R_e\Omega_e^\wedge$:
\begin{equation*}
 \dot{V}_{1}(R_{e})=\left\langle\left\langle \psi(R_{e})\left(R_{e}-R^{T}_{e}\right)^{\vee}, \;-\Omega_{u}(R_{e})\right\rangle\right\rangle <0 ,   
\end{equation*} 
for all $R_{e}\in SO(3)\backslash\mathcal{O}_{R}$ and $\dot{V}_{1}(R_{e})\leq - y_{1}(R_{e})V_{1}(R_{e})$ for all $R_{e}\in \mathcal{U}_{R}$, where 
\begin{equation*}
    y_{1}(R_{e}) \triangleq \frac{\phi(R_{e})}{4\psi(R_{e})\sin\phi(R_{e})} >0,\; \forall R_{e}\in \mathcal{U}_{R}. 
\end{equation*}
This proves that the kinematic control \eqref{eq:KinCtrlLawSO3} also holds Definition \ref{def:KctrLaw}(\ref{Pty3_DefKLaw})-(\ref{Pty4_DefKLaw}). 

Therefore, based on the kinematic control law \eqref{eq:KinCtrlLawSO3}  the following group operation is defined
\begin{align} 
    &r_{1}\star r_{2} \label{eq:GroupOperSO3}  \\
    =&\bigg( R_{1}R_{2},\; \Omega_{1} + \Omega_{2} + \lambda\Omega_{u}(R_{1}) + \lambda\Omega_{u}(R_{2}) -\lambda\Omega_{u}(R_{1}R_{2})  \Big), \notag 
\end{align}
for any $r_{1}=(R_{1},\Omega_{1})$, $r_{2}=(R_{2},\Omega_{2})\in TSO(3)$. Thus, the tangent bundle $TSO(3)\simeq SO(3)\times \mathbb{R}^{3}$ is endowed with a Lie group structure with identity $(I_{3},0_{3\times 1})\in TSO(3)$ and inverse $r^{-1}=(R^{T},-\Omega)\in TSO(3)$,  $\forall r=(R,\Omega)\in TSO(3)$. 
Likewise, given $r_{e}(t)=\left(R_{e}(t),\Omega_{e}(t)\right)\in TSO(3)$ and in view of \eqref{eq:SurfaceErr}, the map $s: TSO(3)\to \mathbb R^3$
\begin{equation}\label{eq:SldsSO3}
    s(r_{e}) = \Omega_{e}(t) + \lambda \Omega_{u}(R_{e}),
\end{equation}
for some scalar $\lambda >0$, defines a Lie subgroup 
\begin{equation}
  H_{R}=\{ r_{e}(t)=\left(R_{e}(t),\Omega_{e}(t)\right)\in TSO(3)\;|\; s(r_{e})=0_{3\times 1} \}, \label{eq:HSO(3)}
\end{equation}
under the group operation \eqref{eq:GroupOperSO3}.

Thus, the tracking controller on $SO(3)$ is obtained from \eqref{eq:TckCtrl} and \eqref{eq:SldsSO3} as
\begin{align}
    \tau_{u} &= \mathbb{J}^{-1}\left( \left(\mathbb{J}\Omega (t)\right)^{\wedge} \left( \lambda\Omega_{u}(R_{e}) -\sigma (t)\right)\right) - \lambda\dot{\Omega}_{u}(R_{e}) +\dot{\sigma}(t) \notag\\
    &\quad -k_{s}s(r_{e}), \label{eq:TckCtrlSO3}
\end{align}
where $k_{s}>0$ is a controller gain. Theorem \ref{thm:TckCtrl} proves that controller \eqref{eq:TckCtrlSO3} in closed loop with the system \eqref{eq:DynSO3} renders the equilibrium point $r_{e}(t)=\left( I_{3},0_{3\times 1}\right)$ almost globally asymptotically stable for all $r_{e}(0)\in SO(3)\backslash \mathcal{O}_{r} \times \mathbb{R}^{3}$, and exponentially stable for all $r_{e}(0)\in \mathcal{U}_{r} \times \mathbb{R}^{3}$.

\begin{rem}
 Note that in applying Theorem \ref{thm:TckCtrl} it should define first a tracking error  using the group operation on the configuration manifold, and then  treat the error dynamics  as a physical system. Otherwise, the sliding surface may not be a Lie subgroup. To see this more clear,  consider $r = \left( R,\Omega \right)$, $r^{-1}_{d}=\left( R^{T}_{r},-\Omega_{r}\right)\in TSO(3)$, then  in the following tracking error may be defined by the group operation \eqref{eq:GroupOperSO3}
\begin{align*}
    r'_{e} & =r^{-1}_{d}\star r = \left( R^{T}_{r}R,\;  - \Omega_{r} + \Omega + \lambda\Omega_{u}(R^{T}_{r}) + \lambda\Omega_{u}(R) \right. \notag\\
    &\quad\quad\quad\quad\quad\quad\quad\quad\quad\left.-\lambda\Omega_{u}(R^{T}_{r}R)  \right) \notag \\
    &= \left( R_{e},\;  - \Omega_{r} + \Omega + \lambda\Omega_{u}(R^{T}_{r}) + \lambda\Omega_{u}(R) -\lambda\Omega_{u}(R_{e})  \right) \notag \\
    &= \left( R_{e},\bar{\Omega}_{e}\right).
\end{align*}
However, 
$H'_{R}=\{ r'_{e}(t)\in TSO(3)\;|\; s(r'_{e})=0_{3\times 1} \}\subset TSO(3)$ 
is not a sliding subgroup for the proposed Morse function $V_{1}(R_{e})$. 
\end{rem}

\subsection{Attitude Tracking on $\mathcal{S}^{3}$}
The set of $\mathbb{R}^{4}$-vectors evolving on the unit sphere $\mathcal{S}^{3}=\left\{ q\in\mathbb{R}^{4} \;|\; q^{T}q = 1\right\}$, with $q=\left[q_{0},\Vec{q}^{T}\right]^{T}\in\mathcal{S}^{3}$, $q_{0}\in [-1,1]$, and $ \Vec{q}\in \mathbb{R}^{3}$, is a Lie group with identity  $\imath = \left[ 1, 0_{1\times 3}\right]^{T}\in\mathcal{S}^{3}$, and inverse $q^{-1}=\left[ q_{0},-\Vec{q}^{T}  \right]^{T}\in\mathcal{S}^{3}$, under the group operation $ (q_{1},q_{2})\mapsto q_{1}\otimes q_{2}\in\mathcal{S}^{3}$
defined as
\begin{equation*}
    q_{1}\otimes q_{2} \triangleq Q(q_{1})q_{2} = \left[ \begin{array}{cc}
      q_{0,1}   & -\Vec{q}^{T}_{1} \\
       \Vec{q}_{1}  & q_{0,1}I_{3} + \Vec{q}_{1}^{\wedge}
    \end{array}\right] \left[ \begin{array}{c}
         q_{0,2}\\
         \Vec{q}_{2}
    \end{array}\right] ,
\end{equation*}
for any $q_{1} = \left[q_{0,1},\Vec{q}^{T}_{1}\right]^{T}$, $q_{2}=\left[q_{0,2},\Vec{q}^{T}_{2}\right]^{T}\in\mathcal{S}^{3}$. The Lie algebra  is  $\mathfrak{s}^{3}=\left\{ \omega \in\mathbb{R}^{4}\;|\; \omega = \left[0,\Omega^{T}\right]^{T}, \Omega\in\mathbb{R}^{3} \right\}$, which  holds $\mathfrak{s}^{3}\simeq \mathbb{R}^{3}$. Its Lie bracket operation corresponds to the cross product in $\mathbb{R}^{3}$. Thus, denote the isomorphism $\overline{\cdot}\vcentcolon \mathbb{R}^{3}\to \mathfrak{s}^{3}$ with the inverse map $\underline{\cdot}\vcentcolon \mathfrak{s}^{3}\to \mathbb{R}^{3}$. 

Rodriguez formula  $q\mapsto R(q) = I_{3}+2q_{0}\Vec{q}^{\wedge} + 2\Vec{q}^{\wedge 2} \in SO(3)$ relates 
each antipodal point $\pm q$ with a physical rotation of a rigid body, 
i.e., $\mathcal{S}^{3}$ double covers the group $SO(3)$. The adjoint action in $\mathcal{S}^{3}$ is defined as $\mathrm{Ad}_{q}\zeta \triangleq q\otimes \Bar{\zeta}\otimes q^{-1} = \overline{R(q)\zeta}$, for any $\zeta\in\mathbb{R}^{3}$.

Given a differentiable curve $q\vcentcolon I \to \mathcal{S}^{3}$ with a left-invariant vector field $\dot{q}(t)\in T_{q}\mathcal{S}^{3}$, and a twice-differentiable reference configuration $q_{r}\vcentcolon I\to \mathcal{S}^{3}$,  $\forall t\in I$, the body angular velocity 
$\overline{\Omega} (t) \triangleq 2q^{-1}(t)\otimes \dot{q}(t) = 2Q^{T}(q(t)) \dot{q}(t)\in\mathfrak{s}^{3}$  
and the reference angular velocity  $\overline{\Omega}_{r}(t)\triangleq 2q^{-1}_{r}(t)\otimes \dot{q}_{r}(t) \in\mathfrak{s}^{3}$ can be defined. 
We consider the following intrinsic tracking error  $ q_{e}(t) \triangleq q^{-1}_{r}(t)\otimes q(t),$
and its left-invariant velocity error
\begin{align}
    \overline{\Omega}_{e}(t) &\triangleq 2q^{-1}_{e}(t)\otimes \dot{q}_{e}(t) = \overline{\Omega}(t) - \overline{\zeta}(t), \label{eq:VelErrS3}\\
    \overline{\zeta}(t) &= \mathrm{Ad}_{q^{-1}_{e}}\Omega_{r}(t),\notag
\end{align}
where $\dot{q}_{e}(t)\in T_{q_{e}}\mathcal{S}^{3}$ is left invariant. 
Propose the Morse function $\mathcal{S}^{3}\ni q \mapsto V_{2}(q)= \frac{1}{\sqrt{2}} \left\| \imath - q\right\| =\sqrt{ 1 -q_{0}}$, which satisfies $V_{2}(q) = 0 \iff q=\imath$, and $V_{2}(q)>0$ $\forall q\in\mathcal{S}^{3}\backslash \{\imath\}$. That is,  function $V_{2}(q)$ has a unique minimum critical zero at identity $\imath\in\mathcal{S}^{3}$ and is strictly positive for any other $q\in\mathcal{S}^{3}$. 
Moreover, it verifies that
\begin{equation*}
    \frac{\mathrm{d}}{\mathrm{d}t}V_{2}(q_{e}) = \frac{-\dot{q}_{0,e}(t)}{2\sqrt{1-q_{0,e}(t)}} = \frac{1}{4\sqrt{1-q_{0,e}(t)}}\Vec{q}^{T}_{e}(t)\Omega_{e}(t),
\end{equation*}
which suggests the following kinematic control law
\begin{equation}\label{eq:KinCtrlLawS3}
    \Omega_{u}(q_{e})\equiv \mathrm{log}(q_{e}) \triangleq \left\{ \begin{array}{ll}
         0_{3\times 1},& q_{e}=\imath ,  \\
         \frac{\arccos (q_{0,e})}{\|\Vec{q}_{e}\|}\Vec{q}_{e} , & q_{e} \neq \imath ,
    \end{array}\right.
\end{equation}
for all $q_{e}(t)\in\mathcal{S}^{3}\backslash \{-\imath\}$. Indeed, when $\Omega_{e}(t)=-\Omega_{u}(q_{e})$, it leads to 
\begin{align*}
    \frac{\mathrm{d}}{\mathrm{d}t}V_{2}(q_{e}) &= -\frac{\arccos (q_{0,e})}{4\sqrt{1-q_{0,e}}}\|\Vec{q}_{e}\| \notag\\
    &= -\frac{\arccos (q_{0,e})}{4\sqrt{1-q_{0,e}}} \sqrt{1-q^{2}_{0,e}} \notag \\
    &= -\frac{\arccos (q_{0,e})}{4\sqrt{1-q_{0,e}}} \sqrt{\left(1+q_{0,e}\right)\left(1-q_{0,e}\right)} \notag\\
    &= -\frac{\arccos (q_{0,e})}{4\sqrt{1-q_{0,e}}} \sqrt{1+q_{0,e}}V_{2}(q_{e}), \notag\\
    &= - y_{2}(q_{e})V_{2}(q_{e}),\notag
\end{align*}
where $y_{2}(q_{e})>0$ for all $q_{e}\in \mathcal{U}_{q}\triangleq \{ q_{e}\in\mathcal{S}^{3}\backslash \{\imath\}\; |\; V_{2}(q_{e})< 2 -\epsilon \}$, for some $\epsilon>0$ arbitrarily small. Consequently, the control law \eqref{eq:KinCtrlLawS3} satisfies all properties of Definition \ref{def:KctrLaw} for the Morse function $V_{2}(q_{e})$.

The kinematic control law \eqref{eq:KinCtrlLawS3} enables the definition of the tangent bundle $T\mathcal{S}^{3} \simeq \mathcal{S}^{3}\times \mathbb{R}^{3}$ as a Lie group under the group operation 
\begin{align} 
 &p_{1}\star p_{2} \label{eq:operation_q} \\
 =& \left( q_{1}\otimes q_{2},\; \Omega_{1} + \Omega_{2} + \lambda\Omega_{u}(q_{1}) + \lambda\Omega_{u}(q_{2}) -\lambda\Omega_{u}(q_{1}\otimes q_{2})  \right), \notag
\end{align}
$\forall p_{1}=(q_{1},\Omega_{1})$, $p_{2}=(q_{2},\Omega_{2})\in T\mathcal{S}^{3}$. Note that 
 the identity is $(\imath,0_{3\times 1})\in T\mathcal{S}^{3}$, and inverse,   $p^{-1}= (q^{-1},-\Omega)\in T\mathcal{S}^{3}$,  $\forall p = (q,\Omega)\in T\mathcal{S}^{3}$. 
Therefore, the map 
\begin{equation}\label{eq:SldSfS3}
    s(p_{e}) = \Omega_{e}(t) + \lambda \Omega_{u}(q_{e}),
\end{equation}
where $p_{e}(t)= (q_{e}(t),\Omega_{e}(t))\in T\mathcal{S}^{3}$, 
defines the sliding Lie subgroup 
\begin{equation}
    H_{q} = \left\{ p_{e}\in T\mathcal{S}^{3}\; |\; s(p_{e})=0_{3\times 1} \right\}. \label{eq:HSO3}
\end{equation}
The attitude tracking controller on $\mathcal{S}^{3}$ is thus defined as \eqref{eq:TckCtrl} using \eqref{eq:VelErrS3}-\eqref{eq:SldSfS3}, which yields 
\begin{align}
    \tau_{u} &= \mathbb{J}^{-1}\left( \left(\mathbb{J}\Omega (t)\right)^{\wedge} \left( \lambda\Omega_{u}(q_{e}) -\zeta (t)\right)\right) - \lambda\dot{\Omega}_{u}(q_{e}) +\dot{\zeta}(t) \notag\\
    &\quad -k_{s}s(p_{e}). \label{eq:TckCtrlS3}
\end{align}
By Theorem \ref{thm:TckCtrl} controller \eqref{eq:TckCtrlS3} in closed loop with system \eqref{eq:ELDynSO3} achieves the asymptotic convergence of $p_{e}(t)\to (\imath,0_{3\times 1})$ for all $p_{e}(0)\in\mathcal{S}^{3}\backslash \{\imath\} \times \mathbb{R}^{3}$, and exponential convergence when $p_{e}(0)\in\mathcal{U}_{q} \times \mathbb{R}^{3}$.

\section{Simulations}\label{sec:simulations}

To illustrate the theoretical results and for comparison, the proposed GSMC \eqref{eq:TckCtrlSO3} was contrasted with two reported controllers: the "linearization"-by-state-feedback-like  (LSF) controller Eq. (26) of \cite{maithripala2006almost}, and the PD+ controller Eq. (23) of \cite{lee2012exponential}.  
For  easy comparison, the applied torque control for each controller is rewritten in terms of  
\begin{equation}\label{eq:si}
    s_{i} = \Tilde{\omega}_{i} + \gamma_{i} \Tilde{\varphi}_{i}, \quad \forall i=1,2,3,
\end{equation}
where $\Tilde{\omega}_{i}\in\mathbb{R}^{3}$ is angular velocity error, $\Tilde{\varphi}_{i}\in\mathbb{R}^{3}$ is attitude error, and $\gamma_{i}>0$ is the control gain.

The proposed GSMC law \eqref{eq:TckCtrlSO3} is expressed as
\begin{align}
    \tau_{1} &= -k_{s}\mathbb{J}s_{1} + F_{1} ,\label{eq:Ctrls1}\\
    s_{1} &= \Omega_{e}(t) + \lambda \Omega_{u}(R_{e}) ,\label{eq:s1}\\
    F_{1} &=  \mathbb{J}\left( -\lambda \dot\Omega_{u}(R_{e}) +\dot{\sigma} \right) + \left(\mathbb{J}\Omega\right)^{\wedge}\left( \lambda\Omega_{u}(R_{e}) -\sigma\right). \notag
\end{align}
Likewise, the LSF controller (26) of \cite{maithripala2006almost} is given by
\begin{align}
    \tau_{2} &= -k\mathbb{J}s_{2} + F_{2} ,\label{eq:Ctrls2}\\
    s_{2} &= \Omega - \Omega_{r} + \frac{\kappa}{k}R^{T}\left( RR^{T}_{r} - R_{r}R^{T} \right)^{\vee} ,\label{eq:s2}\\
    F_{2} &=  \mathbb{J}\dot{\Omega}_{r} - \left(\mathbb{J}\Omega\right)^{\wedge}\left( \Omega\right) - \Omega^{\wedge}\Omega_{r} , \notag
\end{align}
where  $\Tilde{I} = I_{3}$ and $K = \kappa I_{3}$, for some $\kappa >0$. Finally, the PD+ controller (23) of \cite{lee2012exponential} is rewritten as
\begin{align}
    \tau_{3} &= -k_{\Omega}s_{3} + F_{3} ,\label{eq:Ctrls3}\\
    s_{3} &= \Omega_{e} + \frac{k_{R}}{k_{\Omega}}\psi(R_{e})\left( R_{e} - R^{T}_{e} \right)^{\vee} ,\label{eq:s3}\\
    F_{3} &= \mathbb{J}R^{T}R_{r}\dot{\Omega}_{r} +  \left( R^{T}R_{r}\Omega_{r}  \right)^{\wedge}\mathbb{J}R^{T}R_{r}\Omega_{r}. \notag
\end{align}

The inertia tensor was given by
\begin{equation*}
    \mathbb{J} = \left[ \begin{array}{ccc}
         3.6046 & -0.0706 & 0.1491 \\
        -0.0706 & 8.6868 & 0.0449 \\
         0.1491 & 0.0449 & 9.3484
    \end{array} \right] ,
\end{equation*}
while the reference trajectory was calculated as $\Omega_{r}(t) = \left( R^{T}_{r}(t)\dot{R}_{r}(t)\right)^{\vee} = \left[ 0,0.1,0 \right]^{T}$ (rad/s). Furthermore, the initial conditions were chosen as $\Omega(0) = \left(1/\left(2\sqrt{14}\right)\right)\left[ 1,2,3\right]^{T}$ (rad/s), $R_{r}(0) = R_{312}(\pi/4, - \pi , \pi/4)$, where the expression $R_{312}(\varphi , \vartheta ,\psi)$ is a rotation matrix described by the sequence 3-1-2 of Euler angles \citep{shuster1993survey}, and the initial attitude was calculated as $R(0) = R_{r}(0)R_{e}(0)$. 

The simulations were carried out under three scenarios according to the distance between $R_{e}(0)$ and the desired equilibrium $I_{3}$, and to the undesired equilibrium $\mathrm{diag}(1,-1,-1)$ measured by the Morse function $\Psi(R_{e}) \triangleq \frac{1}{2}\mathrm{tr}(I_{3}-R_{e})$ used in \cite{maithripala2006almost}. Therefore, the initial attitudes $R_{e}(0) = R_{312}(0,-0.428\pi,0)$, $R_{e}(0) = R_{312}(0,-0.01\pi,0)\approx I_{3}$, and $R_{e}(0)= R_{312}(0,-0.99\pi,0)\approx \mathrm{diag}(1,-1,-1)$ were assigned. 

Finally, the design parameters for each controller were tuned in such a way that the  energy-consumption level  measured by  $\sqrt{\int^{t}_{0} \tau^{T}_{i}(t)\tau_{i}(t) \mathrm{d}t}$ in the first scenario is the same. The resulting controller gains were $k_{s}=1$, $\lambda = 0.5$ for \eqref{eq:Ctrls1}, $k=1$, $\kappa = 0.5$ for \eqref{eq:Ctrls2}, and $k_{\Omega} = 18.5$, $k_{R} = 9.25$ for \eqref{eq:Ctrls3}. With these design parameters the control gain of \eqref{eq:si} was  $\gamma_{i} = 0.5$ for all $i = 1,2,3$.

\subsection{Scenario 1. Intermediate case.}

\begin {figure}
    \begin{center}
        \includegraphics[scale=0.34,trim = 5mm 5mm 5mm 0mm]{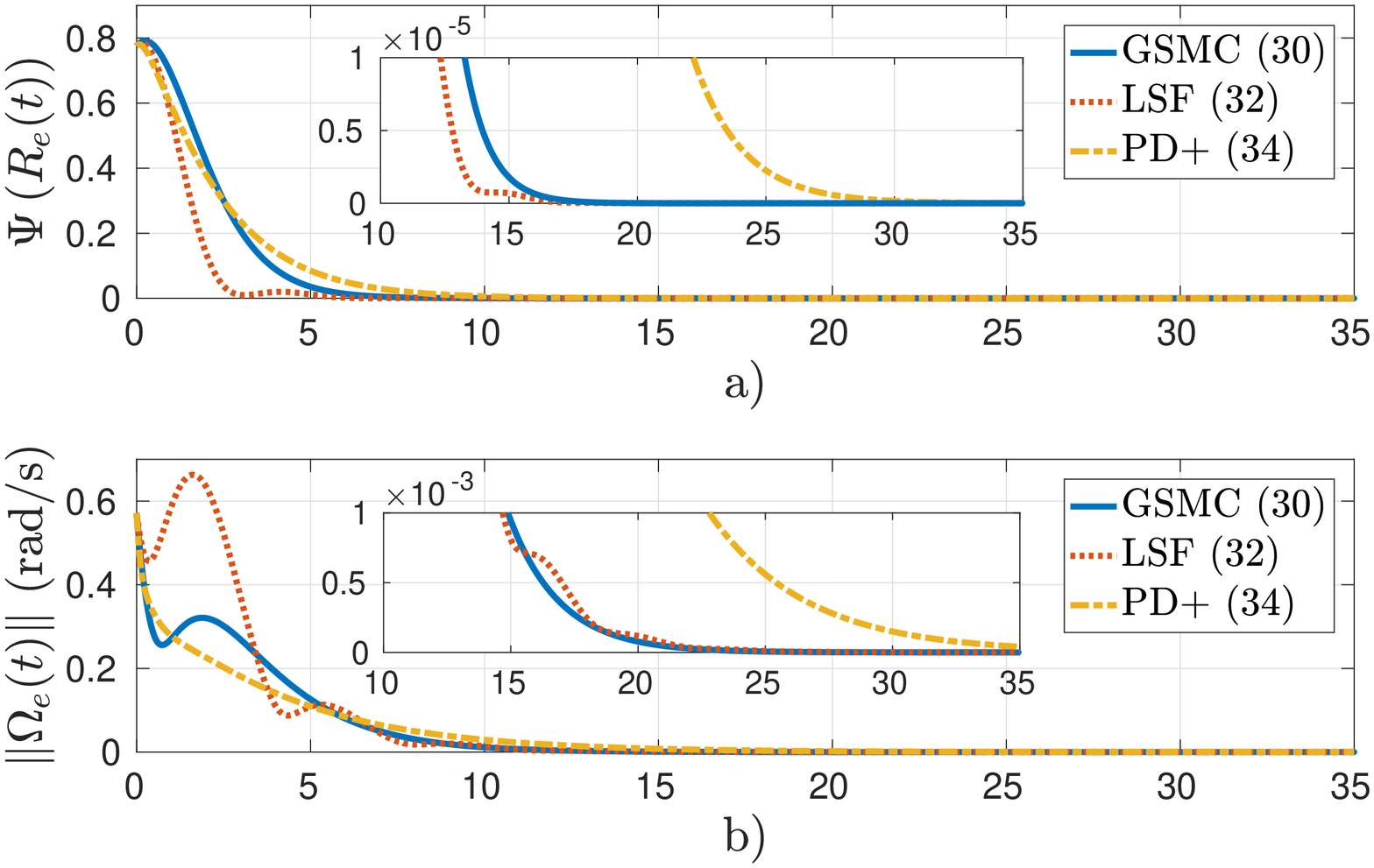}
        \includegraphics[scale=0.34,trim = 5mm 0mm 5mm 0mm]{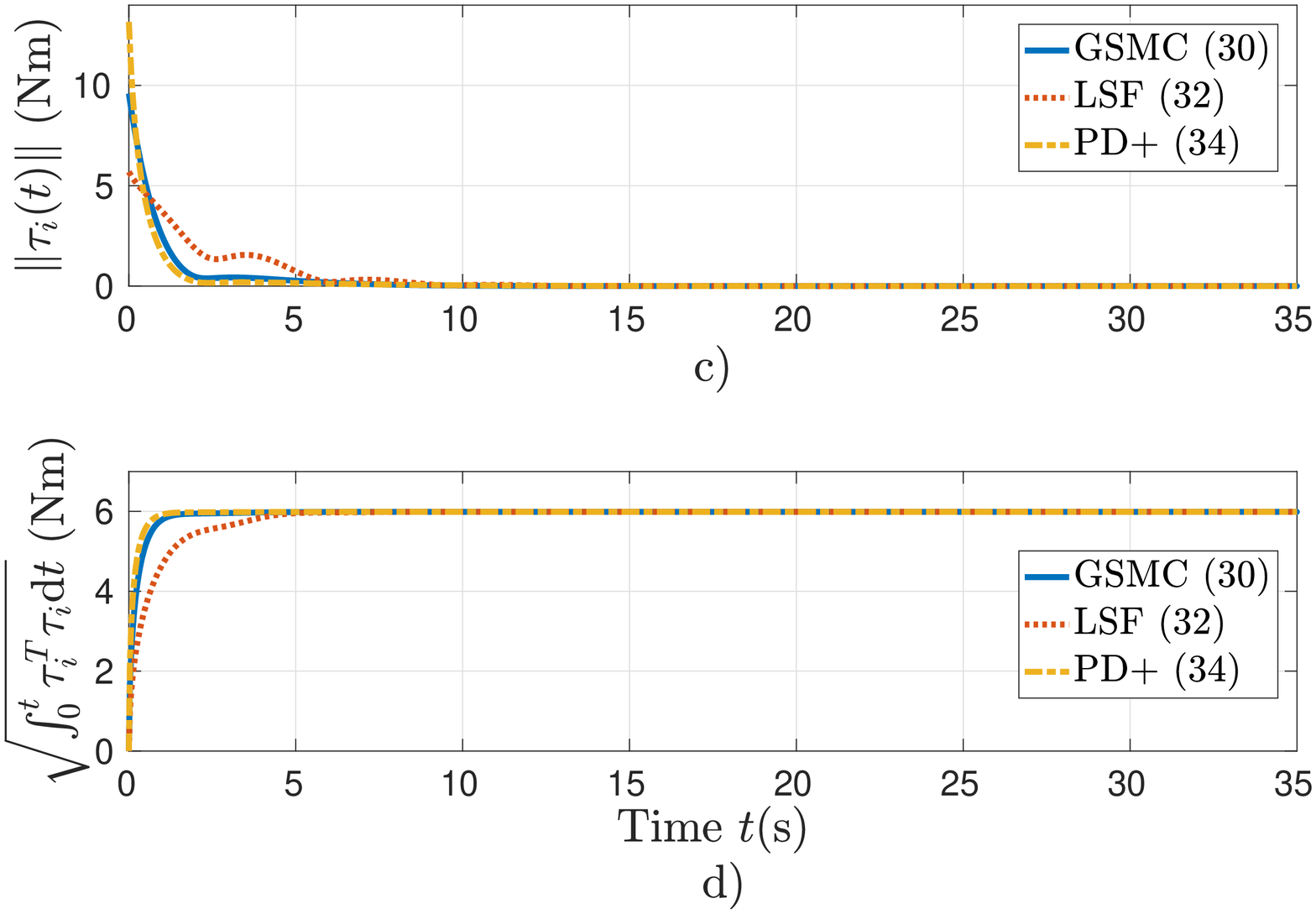}
        \caption{Scenario 1:  Behavior of controllers \eqref{eq:Ctrls1}, \eqref{eq:Ctrls2}, and \eqref{eq:Ctrls3} when the initial attitude error is $R_{e}(0)=R_{312}(0,-0.428\pi,0)$.}
        \label{fig:428Perf}
    \end{center}
\end{figure}

\begin {figure}
    \begin{center}
        \includegraphics[scale=0.34,trim = 5mm 0mm 5mm 0mm]{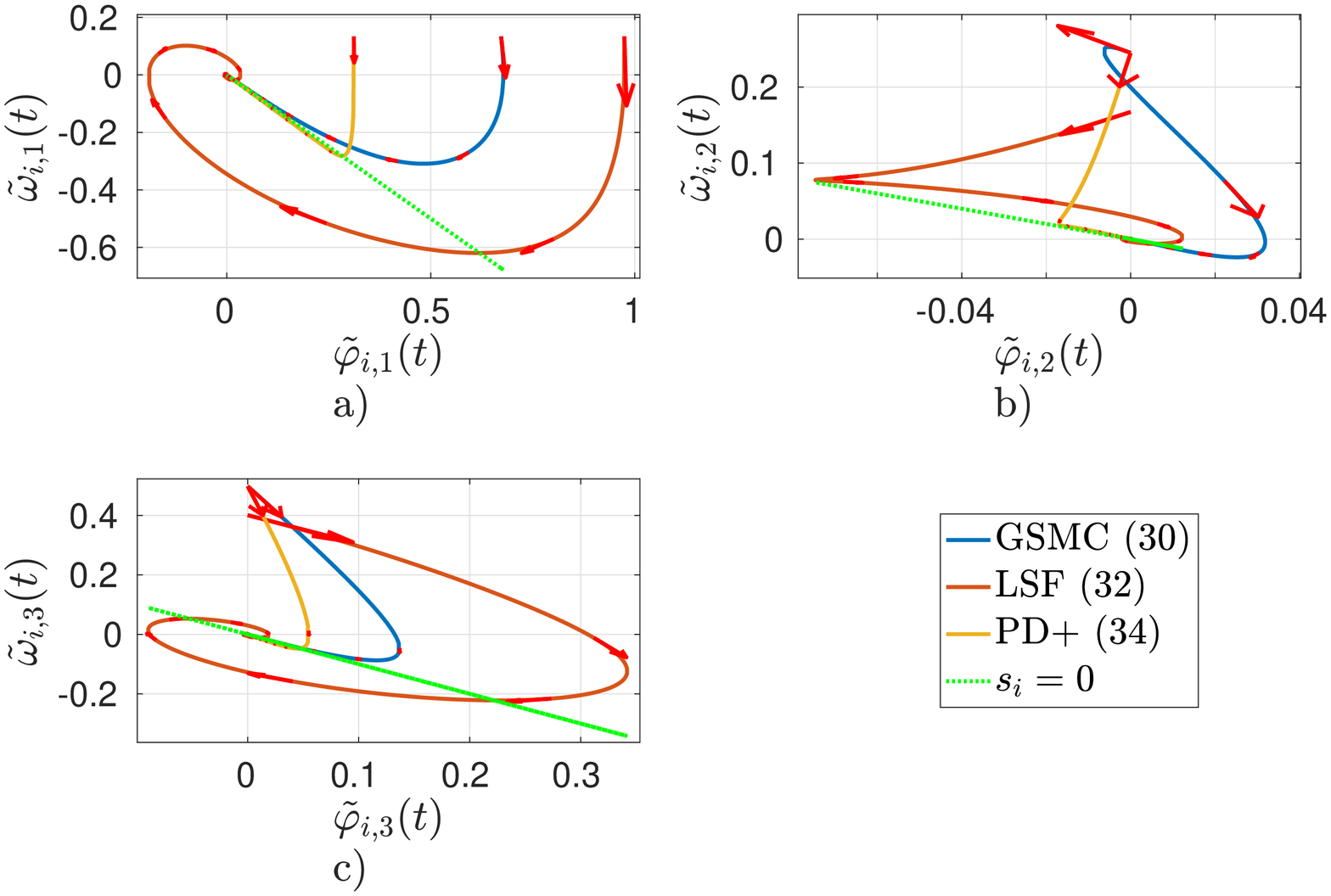}
        \caption{Scenario 1: Behavior of the sliding variable \eqref{eq:s1}, \eqref{eq:s2}, and \eqref{eq:s3} when $R_{e}(0)=R_{312}(0,-0.428\pi,0)$.}
        \label{fig:428SldS}
    \end{center}
\end{figure}

Figure \ref{fig:428Perf} shows the performance of the controllers \eqref{eq:Ctrls1}, \eqref{eq:Ctrls2}, and \eqref{eq:Ctrls3} under the initial condition $R_{e}(0)=R_{312}(0,-0.428\pi,0)$. Fig. \ref{fig:428Perf}(a) shows the attitude error 
$\Psi(R_{e}) \triangleq \frac{1}{2}\mathrm{tr}(I_{3}-R_{e})$, 
it is observed that the proposed controller \eqref{eq:Ctrls1} and controller \eqref{eq:Ctrls2} achieve the convergence $R_{e}\to I_{3}$ in $17$ (s), while controller \eqref{eq:Ctrls3} achieves it in $30$ (s). Fig. \ref{fig:428Perf}(b) illustrates the norm $\|\Omega_{e}(t)\|$ for each controller, where the angular velocity error $\Omega_{e}(t)$ is calculated as \eqref{eq:ErrKynSO3}, it can be seen that controller \eqref{eq:Ctrls3} takes $10$ (s) longer than the other controllers to reach $\Omega_{e}(t)\to 0_{3\times 1}$. Furthermore, Figs. \ref{fig:428Perf}(c) and (d) draw the control effort and the energy consumption respectively, it is observed that, with the selected controller gains, all controllers consume the same amount of energy. Finally, Fig. \ref{fig:428SldS} shows the behavior of the sliding variables \eqref{eq:s1}, \eqref{eq:s2}, and \eqref{eq:s3} compared to $s_{i}=0$ according to \eqref{eq:si}. It is observed that the proposed controller \eqref{eq:Ctrls1} allows convergence $s_{1}\to 0_{3\times 1}$ to complete the reach phase, while the LSF control scheme \eqref{eq:Ctrls2} presents an oscillatory behavior around the equilibrium point, in addition to the PD + controller \eqref{eq:Ctrls3} that follows closely $s_{i}=0$ until it reaches equilibrium.

\subsection{Scenario 2. Starting close to the desired equilibrium point $I_{3}$.}

\begin {figure}
    \begin{center}
        \includegraphics[scale=0.34,trim = 5mm 5mm 5mm 0mm]{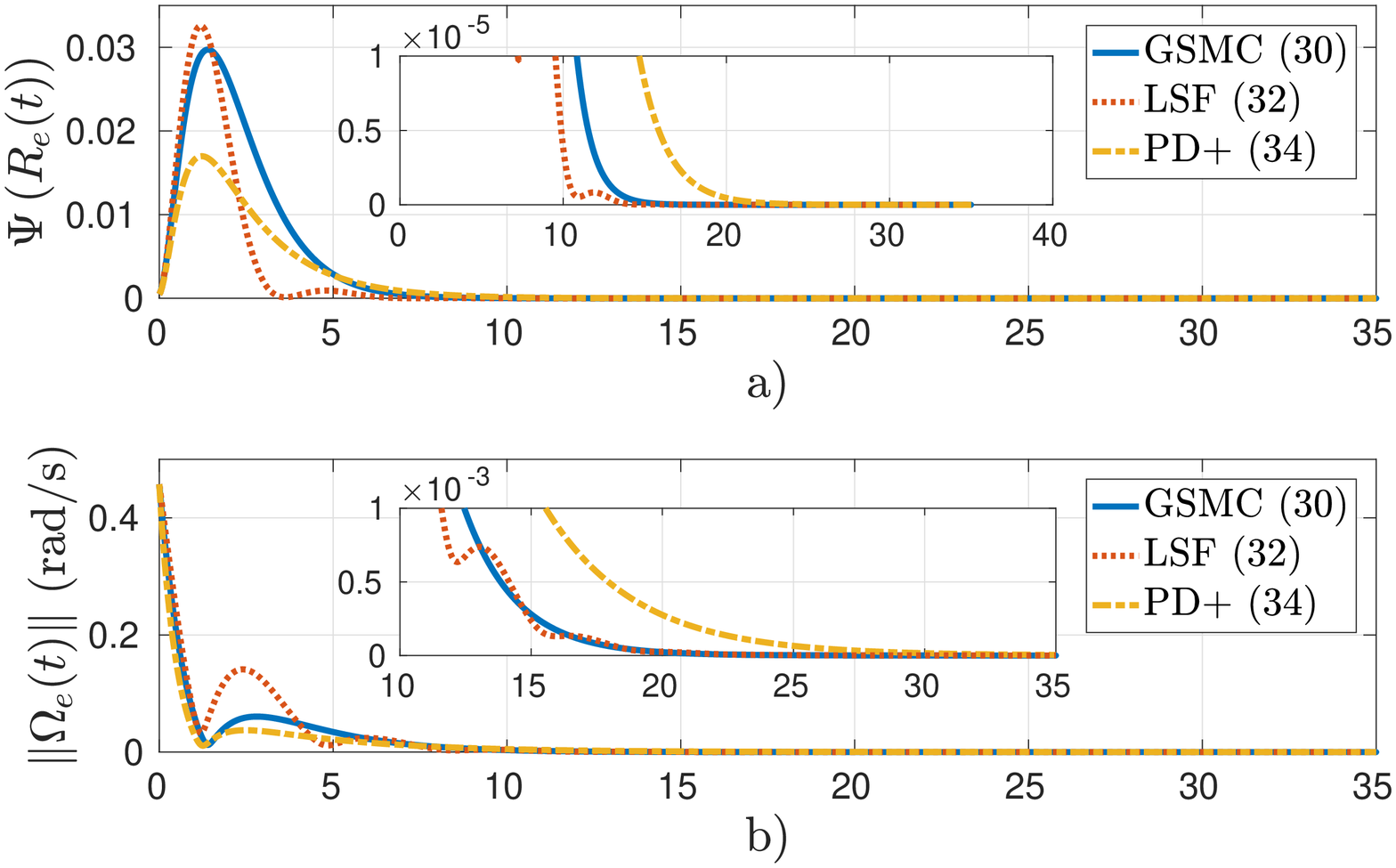}
        \includegraphics[scale=0.34,trim = 5mm 0mm 5mm 0mm]{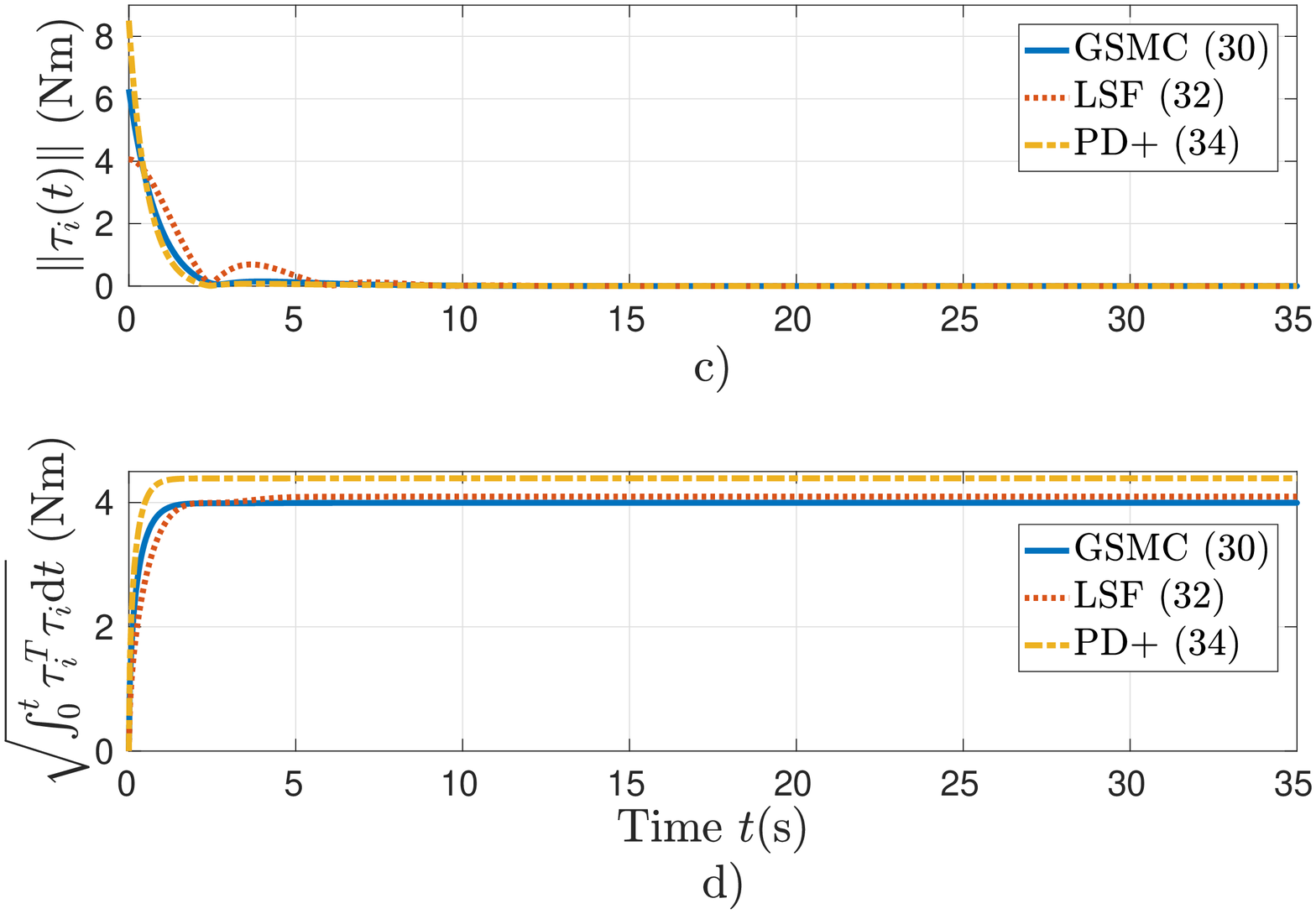}
        \caption{Scenario 2:  Behavior of controllers \eqref{eq:Ctrls1}, \eqref{eq:Ctrls2}, and \eqref{eq:Ctrls3} when the initial attitude error is close to $I_{3}$, i.e.,  $R_{e}(0)=R_{312}(0,-0.01\pi,0)$.}
        \label{fig:001Perf}
    \end{center}
\end{figure}

For this scenario, the initial condition was set to $R_{e}(0)=R_{312}(0,-0.01\pi,0)$, which corresponds to an initial condition close to the desired equilibrium $I_{3}$. Figs. \ref{fig:001Perf}(a) and (b) show that the controllers \eqref{eq:Ctrls1} and \eqref{eq:Ctrls2} reach the desired equilibrium $(R_{e},\Omega_{e}) = (I_{3},0_{3\times 1})$ at the same time $20$ (s), while the controller \eqref{eq:Ctrls3} takes $5$ (s) longer, which coincides with the previous scenario. However, as illustrated in Fig. \ref{fig:001Perf}(d), the proposed controller uses less energy than others to reach the desired equilibrium when the system starts close to the desired equilibrium.

\subsection{Scenario 3. Starting close to the undesired equilibrium point $\mathrm{diag}(1,-1,-1)$.}

\begin {figure}
    \begin{center}
        \includegraphics[scale=0.34,trim = 5mm 5mm 5mm 0mm]{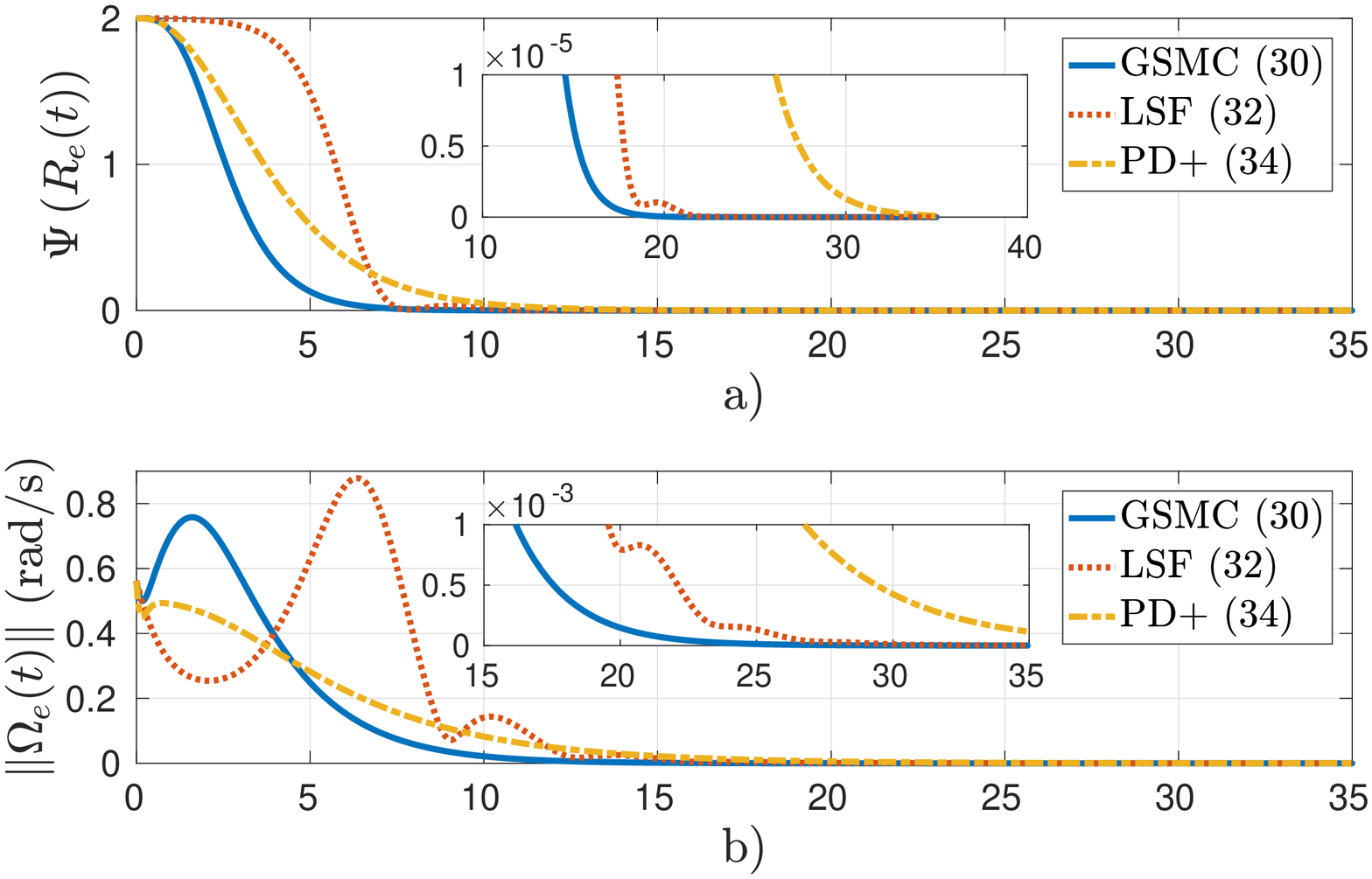}
        \includegraphics[scale=0.34,trim = 5mm 0mm 5mm 0mm]{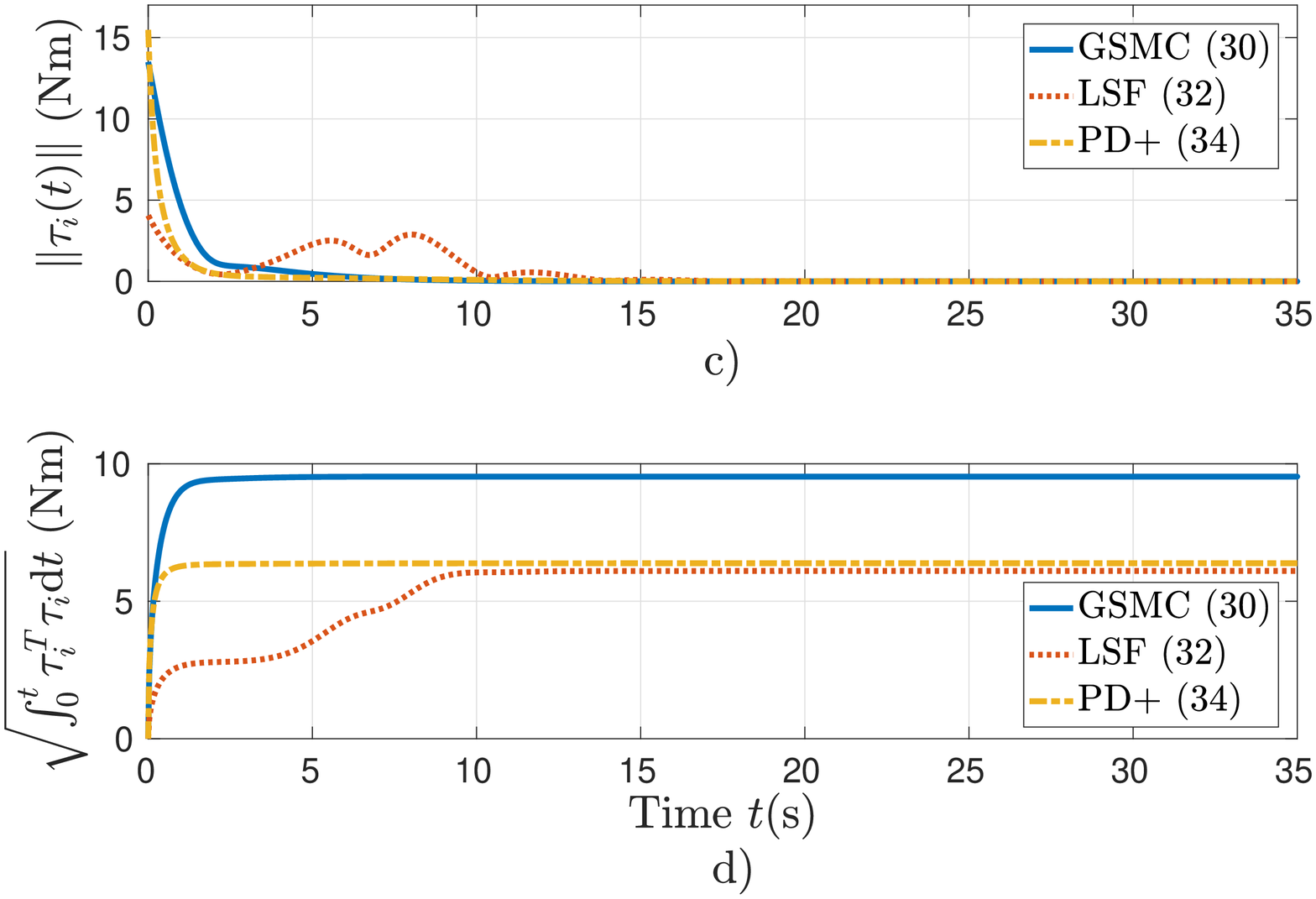}
        \caption{Scenario 3:  Behavior of controllers \eqref{eq:Ctrls1}, \eqref{eq:Ctrls2}, and \eqref{eq:Ctrls3} when the initial attitude error is close to $\mathrm{diag}(1,-1,-1)$), i.e.,  $R_{e}(0)=R_{312}(0,-0.99\pi,0)$.}
        \label{fig:099Perf}
    \end{center}
\end{figure}

Figure \ref{fig:099Perf} displays the performance of the controllers starting close to the undesired equilibrium point $\mathrm{diag}(1,-1,-1)$, i.e., $R_{e}(0)=R_{312}(0,-0.99\pi,0)$. It is observed in Fig. \ref{fig:099Perf}(a) that the proposed controller \eqref{eq:Ctrls1} and the PD+ controller \eqref{eq:Ctrls1} present a delay of $1$ (s) before beginning the convergence of $R_{e}\to I_{3}$, however, the LSF controller \eqref{eq:Ctrls2} has the longest delay of $2.5$ (s). Notice that the proposed control scheme allows a faster convergence to the desired equilibrium point (Figs. \ref{fig:099Perf}(a) and (b)) at a cost of more energy consumption (Figs. \ref{fig:099Perf}(c) and (d)).

\section{Conclusions} \label{sec:Conclusions}

This paper presented a geometric sliding mode control for fully actuated mechanical systems evolving on Lie groups, generalizing the conventional sliding mode control in Euclidean spaces. 
It was shown that the sliding surface (a Lie subgroup) is immersed in the state space (a Lie group)  of the system dynamics, and the tracking is achieved by first driving the trajectories of the system to the sliding subgroup and then converging to the group identity of the reduced dynamics restricted on the sliding subgroup, like sliding mode control designs for systems evolving on Euclidean spaces. An application of the result to attitude control was presented for the rotation group $SO(3)$ and the unit sphere $\mathcal S^3$. The simulation results illustrated the scheme and compared it with similar control designs in the literature.  



\bibliographystyle{agsm}
\bibliography{GSMC-v1}           



\appendix
\section{A passivity-like lemma} \label{AppA}   
\begin{lem} \label{lm:skysym}
Given that the inner product on $\mathfrak{g}$ is a symmetric bilinear map, for any $\zeta,\eta \in \mathfrak{g}$, it holds 
\begin{align*}
    \mathbb{I}\left( \zeta,\; \mathbb{I}^{\sharp}\left( \mathrm{ad}^{*}_{\zeta}\mathbb{I}^{\flat}(\eta)\right)\right) &= \langle \mathbb{I}^{\flat}\left( \mathbb{I}^{\sharp}\left( \mathrm{ad}^{*}_{\zeta}\mathbb{I}^{\flat}(\eta)\right)\right) ; \zeta \rangle \\
    &= \langle  \mathrm{ad}^{*}_{\zeta}\mathbb{I}^{\flat}(\eta) ; \zeta \rangle \\
    &= \langle \mathbb{I}^{\flat}(\eta) ; \left[ \zeta ,\zeta \right] \rangle \\
    &= 0,
\end{align*}
because of the skew symmetry of the Lie bracket operation $[\cdot,\cdot]\in\mathfrak{g}$.

\end{lem}
\end{document}